\documentclass[11pt]{amsart}

\RequirePackage{amsmath, amssymb, amsfonts, amsthm}

\usepackage{eucal} 

\RequirePackage{enumerate}  

\newcounter{pulse}
\numberwithin{pulse}{section}
\theoremstyle{plain}
\newtheorem{thm}[pulse]{Theorem}
\newtheorem{propn}[pulse]{Proposition}
\newtheorem{lemma}[pulse]{Lemma}

\theoremstyle{definition}
\newtheorem{defn}[pulse]{Definition}
\newtheorem*{defn*}{Definition}
\newtheorem*{notn*}{Notation}
\newtheorem{eg}[pulse]{Example}
\theoremstyle{remark}
\newtheorem*{rem}{Remark}

\newcommand{\mc}{\mathcal}
\newcommand{\mf}{\mathfrak}

\newcommand{\al}{\alpha}
\newcommand{\lm}{\lambda}
\newcommand{\sig}{\sigma}

\RequirePackage[PS,nohug]{diagrams}

\newcommand{\ackndiag}{Uses Paul Taylor's {\tt diagrams.sty} macros.}

\newcommand{\wtild}{\widetilde}
\newcommand{\blob}{\bullet}
\newcommand{\blank}{\underline{\quad}}
\newcommand{\iso}{\cong}
\newcommand{\Ind}{{\mathbb I}} 

\newcommand{\Ker}{\mathop{\rm ker}}
\newcommand{\Image}{\mathop{\rm im}}

\newcommand{\diag}{{\rm diag}} 

\newcommand{\pcat}[1]{{\sf {#1}}} 
\newcommand{\vctG}[1]{{\boldsymbol{#1}}} 
\newcommand{\vctR}[1]{{\bf {#1}}} 

\newcommand{\st}{\;:\;}
\newcommand{\defeq}{:=}
\newcommand{\dt}[1]{{\it #1}\/}  
\renewcommand{\emph}[1]{{\sl #1\/}}  
\newcommand{\sid}{{\sf id}}  
\newenvironment{inline}[1]
{\vspace{0.5em} \noindent {\bf {#1}} }{\vspace{0.5em}} 

\newcommand{\BAlg}{\pcat{BAlg}} 

\newcommand{\norm}[1]{\Vert{#1}\Vert}
\newcommand{\Norm}[1]{\left\Vert{#1}\right\Vert}

\newcommand{\tp}{\mathop{\scriptstyle\otimes}}
\newcommand{\ptp}{\mathop{\scriptstyle\widehat{\otimes}}}
\newcommand{\ptpR}[1]{\mathop{\scriptstyle\widehat{\otimes}}\limits_{#1}}

\newcommand{\Cplx}{{\mathbb C}}
\newcommand{\Nat}{{\mathbb N}}
\newcommand{\Z}{{\mathbb Z}}

\newcommand{\lp}[2][]{\ell^{#2}_{#1}}
\newcommand{\lpsum}[1]{\overset{(#1)}{\bigoplus}}

\newcommand{\id}[1][]{{\sf 1}_{#1}} 
\newcommand{\fu}[1]{{#1}^\#}    

\newcommand{\B}{{\mathcal B}}
\newcommand{\X}{{\mf X}}

\newcommand{\Ho}[3][]{{\mathcal {#2}_{#3}^{#1}}}  
\newcommand{\Co}[3][]{{\mathcal {#2}^{#3}_{#1}}} 

\newcommand{\bdy}{{\sf d}}
\newcommand{\face}[2][]{\partial^{#1}_{#2}} 

\newarrow{Sub}       C--->  
\newarrow{Eq} ===== 
\newarrow{Dash} {}{dash}{}{dash}>

\newcommand{\LR}[2]{\pile{\lTo^{#1} \\ \rTo_{#2} }}

\newcommand{\SA}[1][]{{\mathfrak {alg}}_{#1}}

\newcommand{\SAcat}{\mathbf{SA}}  

\newcommand{\Base}{{\rm base}}

\newcommand{\inc}[2][]{\iota^{#1}_{#2}}  
\newcommand{\Vinc}[2][]{\vctG{\iota}^{#1}_{\vctR{#2}}} 

\newcommand{\ev}[1]{\widehat{\vctR{#1}}} 
\newcommand{\Tran}[1]{\mathop{\sf Tr}\nolimits^{#1}} 

\newcommand{\Cnd}[1]{{\sf (#1)}}  

\newcommand{\Ax}{A\ev{x}} 

\title[Strong semi\-lattices of Banach algebras]{Simplicial homology of strong semi\-lattices of Banach algebras}
\author{Yemon Choi}
\address{School of Mathematics and Statistics, University of Newcastle, Newcastle upon Tyne, NE1 7RU, England}
\curraddr{Department of Mathematics, University of Manitoba, Winnipeg, Manitoba, Canada R3T 2N2}
\email{y.choi.97@cantab.net}
\thanks{\ackndiag}

\begin{document}
\maketitle
\begin{abstract}
Certain semi\-groups are known to admit a `strong semi\-lattice decomposition' into simpler pieces. We introduce a class of Ban\-ach algebras that generalise the $\lp{1}$-convolution algebras of such semi\-groups, and obtain a dis\-inte\-gration theorem for their simp\-licial homo\-logy.
Using this we show that for any Clifford semigroup $S$ of amenable groups, $\lp{1}(S)$ is simplicially trivial: this generalises results in~\cite{YC_GMJ}. Some other applications are presented.
\end{abstract}

\begin{section}{Introduction}
The connections between structural properties of a semi\-group and cohomo\-logical properties of its $\lp{1}$-convolution algebra are still not fully understood. However, much work has been done (and much is known) for the case of group algebras. In attempting to extend existing techniques and results to more general semi\-groups, a promising class to consider first is that of so-called \dt{Clifford semi\-groups}. Partial results have been obtained in low dimensions \cite{BowDunc,GouPouWh} but hitherto there seems to have been no systematic attack on the higher-degree cohomo\-logy problems.

Clifford semi\-groups are known to arise as \dt{strong semi\-lattices of groups} (terminology which will be explained below). Recently it was established in \cite{YC_GMJ} that the $\lp{1}$-convolution algebras of semi\-lattices are simplicially trivial; in this article we generalise the techniques of \cite{YC_GMJ} to deal with Clifford semi\-groups.
More precisely, for semi\-groups $S$ which admit a `strong semi\-lattice decomposition', we obtain a `disintegration' theorem for the simplicial homology of $\lp{1}(S)$ (see Theorems \ref{t:mainresult}, \ref{t:mainthm} below). As a corollary we prove that $\lp{1}(S)$ is simplicially trivial whenever $S$ is a Clifford semi\-group with amenable constituent groups. (If $S$ is a commutative Clifford semi\-group, then as in \cite{YC_GMJ} we may go on to deduce that cohomo\-logy of $\lp{1}(S)$ vanishes in degrees 1 and above for any symmetric coefficient bimodule.)

Our disintegration result also covers other instances of semi\-groups with a strong semi\-lattice decomposition: in particular we show that the $\lp{1}$-convolution algebras of normal bands are simplicially trivial.

So far we have not mentioned the `strong semi\-lattices of Banach algebras' of the title. These are most succinctly described, in the language of category theory, as semi\-lattice-shaped diagrams in the category of Banach algebras and contractive homo\-morphisms. We shall give a more concrete description below; for now it suffices to note that this construction
\begin{enumerate}[(i)]
\item gives rise to Banach algebras that include $\lp{1}$-convolution algebras of Clifford semi\-groups as special cases, and
\item seems to provide the correct level of generality for our disintegration theorem \ref{t:mainthm}.
\end{enumerate}
Our definition is slightly more general than the original one given in \cite[\S5.1]{YC_PhD}.

We remark that in a recent preprint \cite{GHS_slatt},
Ghandehari, Hatami and Spronk have independently introduced a more general notion of a Banach algebra 
\dt{graded over a semi\-lattice}\/: however, they restrict attention to the case where the semilattice is finite (since they are primarily interested in amenability questions), while the present paper is motivated almost exclusively by the case where the indexing semilattice is infinite.

Before starting on the definitions, lemmas and theorems some general remarks on this article are in order. We will need to draw on background material from functional analysis: in particular we assume the reader is familiar with the projective tensor product of Banach spaces and elementary Banach space properties of $\lp{1}$\/. For homology and cohomo\-logy of Banach algebras we refer the reader to the two classic sources \cite{Hel_HBTA,BEJ_CIBA}.

We also assume familiarity with some basic category theory, such as may be found in \cite{Mac_CWM}; this is not vital for the proofs to follow but provides some context for the underlying strategy.

\end{section}

\begin{section}{Definitions and preliminaries}
Since our eventual interest is in $\lp{1}$-convolution algebras of semi\-groups, we start by recalling some basic notions from semi\-group theory.
\begin{defn}
A \dt{semi\-lattice} is a commutative semi\-group where each element is idempotent.

Each semi\-lattice has a canonical partial order on it, defined by
\[ f\preceq e \iff ef = f \]
and so may be regarded as a poset in a natural way.
\end{defn}

The following example will play an important role later; we introduce it here for emphasis. (Note that it differs slightly from the corresponding example in \cite{YC_GMJ}.)
\begin{eg}[Free semi\-lattices]\label{eg:ff-slatt}
Let $\Ind$ be a finite set and consider $2^\Ind\setminus\{\emptyset\}$, the collection of all non-empty subsets of $\Ind$\/. We may equip $2^\Ind\setminus\{\emptyset\}$ with the structure of a semi\-lattice, by taking the product of two members $A$ and $B$ to be $A\cup B$; note that the canonical partial order is then given by
\[ C \preceq D \iff C \cup D = C \iff C\supseteq D \]
We refer to the semi\-lattice thus obtained as the \dt{free semi\-lattice generated by $\Ind$}. 
(One could extend this definition to the case where $\Ind$ is infinite, but for our purposes the finite case is enough.)
\end{eg}

Semilattices have been used in the structure theory of semi\-groups; one fundamental instance of this is the concept of a strong semi\-lattice of semi\-groups, which we now define.

Let ${\mc S}$ be a non-empty class of semi\-groups which is closed under taking semi\-group homo\-morphisms. Alternatively, ${\mc S}$ may be viewed as a full subcategory of the category $\pcat{SGp}$ of semi\-groups and semi\-group homo\-morphisms.
\begin{defn}[cf.~{\cite[Ch.~IV]{How_fund-sgp}}]
A \dt{semi\-lattice diagram in ${\mc S}$} consists of: a semi\-lattice $L$; a family $(S_e)_{e\in L}$ of semi\-groups in ${\mc S}$; and semi\-group homo\-morphisms $\varphi_{f,e}:S_e \to S_f$ for each pair $(e,f) \in L\times L$ such that $e\succeq f$, which satisfy the following compatibility conditions:
\begin{itemize}
\item $\varphi_{e,e}$ is the identity homo\-morphism on $S_e$\,;
\item if $e \succeq f \succeq g$ in $L$ then $\varphi_{g,f}\circ\varphi_{f,e}=\varphi_{g,e}$\/.
\end{itemize}
\end{defn}

To each such semi\-lattice diagram we may associate a semi\-group as follows: the underlying set is defined to be the disjoint union $\coprod_{e\in L} S_e$\,; and multiplication is defined by
\[ x\cdot y \defeq \varphi_{ef,e}(x)\varphi_{ef,f}(y)\qquad(x\in S_e, y \in S_f) \]
or, more pictorially, by
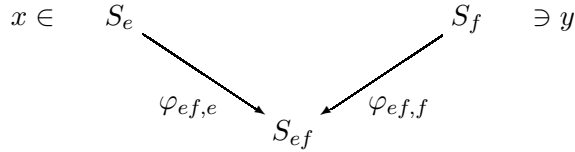
\[ \begin{diagram}[tight,height=2em]
x\in & S_e & & & & S_f & \ni y \\
& & \rdTo_{\varphi_{ef,e}} & & \ldTo_{\varphi_{ef,f}}&  \\
& &  & S_{ef} & & \\
\end{diagram} \]
Any semi\-group $T$ of this form is said to be a \dt{strong semi\-lattice of semi\-groups} (of \dt{type} ${\mc S}$ and \dt{shape} $L$). In this case we refer to $L, (S_e), (\varphi_{f,e})$ as the \dt{decomposition data} for $T$\/.

\begin{eg}
Let ${\mc G}$ be the class of groups. Then a strong semi\-lattice of semi\-groups of type ${\mc G}$, or more concisely a \dt{strong semi\-lattice of groups}, is called a \dt{Clifford semi\-group}. (One may also define Clifford semi\-groups in intrinsic terms, in which case the existence of suitable decomposition data for a given Clifford semi\-group is a theorem rather than a tautology. This can be found in most reference texts on semi\-groups, e.g.~\cite[\S IV.2]{How_fund-sgp}.
\end{eg}

If $T$ is a strong semi\-lattice of semi\-groups, with decomposition data $L$, $(S_e)$ etc.\ then clearly each $\lp{1}(S_e)$ is a subalgebra of $\lp{1}(T)$, and in a loose sense the Banach algebra $\lp{1}(T)$ is built out of the Banach algebras $\lp{1}(S_e)$\/. This observation motivates the following definition, which is a mild generalisation of that in \cite[Ch.~5]{YC_PhD}.

The objects of study in this article are `semi\-lattice-shaped diagrams' in the category $\BAlg_1$ of Banach algebras and contractive homo\-morphisms. An algebraist's definition (following the language of \cite{Mac_CWM}, say) might be as follows:

\begin{inline}{Provisional definition.}
A \dt{strong semi\-lattice of Banach algebras}, or a \dt{semi\-lattice diagram in $\BAlg_1$}, is a pair $(L,A)$ where $L$ is a semi\-lattice and $A$ is a functor from the small category $(L,\succeq)$ to the category $\BAlg_1$\/.
\end{inline}

This definition relies on the interpretation of a poset as a special kind of small category: see \cite[\S 1.2]{Mac_CWM}, for example. The
`proper' definition that follows is merely an unpacking of the terminology in the `provisional' definition just given.
\begin{defn}
A \dt{semi\-lattice diagram in $\BAlg_1$} consists of: a semi\-lattice~$L$; a family $(A_e)_{e \in L}$ of Banach algebras; and \emph{contractive} algebra homo\-morphisms $\phi_{f,e}:A_e \to A_f$ for each pair $(e,f) \in L\times L$ such that $e\succeq f$, which satisfy the following compatibility conditions:
\begin{itemize}
\item $\phi_{e,e}$ is the identity homo\-morphism on $A_e$\,;
\item if $e \succeq f \succeq g$ in $L$ then $\phi_{g,f}\circ\phi_{f,e}=\phi_{g,e}$\/.
\end{itemize}

We shall use the following notational shorthand: the expression `let $(L,A)\in \SAcat$' is henceforth used as an abbreviation for the phrase `let $(L,A)$ be a semi\-lattice diagram in $\BAlg_1$', although on occasion we shall revert to the longer phrase in order to state certain results.
\end{defn}

\begin{rem}
Note that in saying `let $(L,A)$ be a semi\-lattice diagram \ldots' we are suppressing explicit mention of the family $(\phi_{f,e})$\/. This abuse of notation should not lead to confusion.
\end{rem}

\begin{rem}
Our terminology perhaps deserves some remarks:
\begin{enumerate}[(a)]
\item One may think of such a pair $(L,A)$ as a `projective system of Banach algebras' indexed by the partially ordered set $(L,\succeq)$\/. However, such terminology does not seem apposite, as we shall consider neither the inductive nor projective limit of such a system, and since we later use properties of $L$ that are not shared by arbitrary posets.
\item The terminology `$(L,A)\in\SAcat$' strongly suggests that we might view strong semi\-lattices of Banach algebras as objects of a category $\SAcat$\/. While this is easy and not at all deep, we feel that such extra abstraction would not make the results of this paper any easier to follow.
\end{enumerate}
\end{rem}

Given $(L,A)\in\SAcat$, we can endow the $\lp{1}$-sum $\lpsum{1}_{e \in L} A_e$ with a multiplication that turns it into a Banach algebra. Before doing so, we shall introduce some basic notation regarding $\lp{1}$-sums: the point of our pedantry will be explained in due course.

\begin{notn*}
Let $\Ind$ be an indexing set and let $(X_i)_{i\in\Ind}$ be a family of Banach spaces. We write $\X$ for the $\lp{1}$-sum $\lpsum{1}_{i\in\Ind} X_i$, and for each $j\in\Ind$ we let $\inc[X]{j}$ denote the canonical inclusion map $X_j \rSub \X$\/.
If $x \in \coprod_{j\in\Ind} X_j $ then we shall sometimes write $\Base(x)$ for the unique $i \in \Ind$ such that $x \in X_i$\/. Elements of the form $\inc[X]{j}(v)$, where $v\in X_j$, will be called \dt{block elements}.
\end{notn*}
It is evident that two bounded linear maps with domain $\X$ coincide if and only if they agree on all block elements; we shall use this repeatedly below without proof.

\begin{defn}
Let $(L,A)\in \SAcat$\/. We can form a Banach algebra $\SA$, whose underlying vector space is the $\lp{1}$-sum $\lpsum{1}_{e \in L} A_e$ and whose multiplication is defined by the following rule:
\[ \left(\sum_{e\in L}\inc[A]{e} a_e \right)\cdot\left(\sum_{f \in L}\inc[A]{f} b_f\right) \defeq \sum_{g \in L} \left[\sum_{(e,f) \in L\times L\st ef=g} \phi_{g,e}(a_e)\phi_{g,f}(b_f) \right]\]
We say that $\SA$ is the \dt{convolution algebra of $(L,A)$}. If we need to make the dependence on $A$ and $L$ explicit we write $\SA[L,A]$\/.
\end{defn}
Note that even if each $A_e$ is unital, $\SA$ need not be.

\medskip\noindent{\bf Simple examples.}
The following special cases of our construction should be kept in mind, and hopefully serve to clarify the general picture.
\begin{enumerate}[(a)]
\item Fix a Banach algebra ${\mf B}$; let $A_e={\mf B}$ for all $e \in L$; and let each transition map $\phi_{f,e}$ be the identity map on ${\mf B}$\/. Then there is an isomorphism $\lp{1}(L)\ptp {\mf B} \to \SA[L,A]$, defined by sending  $e\ptp b$ to $\inc{e}b$ for every $e \in L$ and $b \in {\mf B}$\/.
\\
\emph{In particular, if ${\mf B}=\Cplx$ then $\SA\iso \lp{1}(L)$, the usual convolution algebra of the semi\-group $L$\/.}
\item At the other extreme, suppose that $L=\{\id\}$, the trivial semi\-lattice containing only one element. Then clearly $\SA[L,A]\iso A_{\id}$\/.
\item Let $L=\{\id, e\}$ be the two-element semi\-lattice consisting of an identity element $\id$ and an idempotent $e$ distinct from $\id$: then $\id \succeq e$\/. Let ${\mf B}$ be any \emph{unital} Banach algebra and define the functor $A:L \to \BAlg_1$ by taking $A_{\id}=\Cplx$, $A_e ={\mf B}$ and $\phi_{e,\id}:\Cplx \to {\mf B}$ to be the homo\-morphism $\lm \mapsto \lm\id[{\mf B}]$\/.
\begin{diagram}[tight,height=2em]
\Cplx & & {\mf B} \\
\id & \rTo & e
\end{diagram}
Then $\SA[L,A]$ is by definition isomorphic as a Banach space to the $\lp{1}$-sum $\Cplx\oplus {\mf B}$; and the multiplication on $\SA[L,A]$ is given by
\[  (\lm, b) \cdot (\mu, c) = (\lm\mu, \lm c+\mu b + bc) \quad\quad(\lm,\mu\in\Cplx; b,c \in {\mf B}) \;.\]
\emph{Thus in this instance $\SA[L,A]$ is nothing but the forced unitisation of ${\mf B}$\/.}
\end{enumerate}

The following example provides the main motivation for our definition of the convolution algebra.
\begin{eg}[$\lp{1}$-algebras of Clifford semi\-groups]
Let $L$ be a semi\-lattice and let $(G_e)_{e\in L}$ be a semi\-lattice diagram of groups. Let $G=\coprod_{e\in L}G_e$ be the corresponding semi\-lattice of groups.

Applying the $\lp{1}$-group algebra functor to each $G_e$, we obtain a semi\-lattice diagram in $\BAlg_1$ with constituent Banach algebras $A_e=\lp{1}(G_e)$\/. It is straightforward to check that $\SA[L,A]$ coincides with the convolution algebra~$\lp{1}(G)$\/.
\end{eg}

\subsection*{$\SA[L,A]$ as an $\lp{1}(L)$-algebra}
The convolution algebras $\SA[L,A]$ have extra structure, whose precise formulation requires a digression into some of the general theory of Banach modules.

\begin{defn}\label{dfn:K-BA}
Recall that if $K$ is a commutative, unital Banach algebra then a \dt{Banach $K$-algebra} is a Banach algebra $A$, equipped with the structure of a symmetric Banach $K$-bimodule, such that for each $a\in A$ both left and right multiplication by $a$ are $K$-module maps.

We say $A$ is a \dt{unit-linked} Banach $K$-algebra if it is unit-linked as a Banach $K$-bimodule. 
\end{defn}

\begin{rem}
Note that a unit-linked, Banach $\Cplx$-algebra is nothing but a Banach algebra in the usual sense.
\end{rem}

It is easily verified that $\SA[L,A]$ has the natural structure of a unit-linked Banach $\lp{1}(L)$-algebra, as follows. We first define a symmetric action of $\lp{1}(L)$ on the vector space $\lpsum{1}_{e \in L} A_e$ by
\[ e\cdot \inc{i}(a) = \inc{i}(a)\cdot e \defeq \inc{ie}\phi_{ie,i}(a) \quad\quad(i,e\in L\/; a \in A_i)\/. \]
Straightforward but tedious computations then show that for any $x\in \SA[L,A]$, the operators of left and right multiplication by $x$ are both $\lp{1}(L)$-bimodule maps.

\begin{rem}
If each $A_e$ is a unital Banach algebra and each $\phi_{f,e}$ a unital homo\-morphism, then there is a canonical, isometric embedding of the vector space $\lp{1}(L)$ into $\lpsum{1}_{e \in L} A_e$, defined by $e \mapsto \iota_e(\id[A_e])$\/. In this case, the way that we have defined multiplication in the convolution algebra $\SA[L,A]$ ensures that this isometric embedding is an algebra homo\-morphism from the convolution algebra $\lp{1}(L)$ \emph{into the centre of $\SA[L,A]$}.
\end{rem}

It is clear that if the semi\-lattice $L$ has an identity $\id$, then $\iota_{\id}(\id)$ is an identity element for the Banach algebra $\SA[L,A]$; if this is the case then the $\lp{1}(L)$-module strcuture on $\SA[L,A]$ is unit-linked. Later on it will be convenient, for technical reasons, to know that a stronger result is true.
\begin{lemma}\label{l:unitalSA}
Let $(L,A) \in \SAcat$\/. Suppose the convolution algebra $\lp{1}(L)$ has an identity element $u=\sum_{e \in L} \lm_e e$; then for every $x \in\SA[L,A]$ we have
\[ u\cdot x = x = x\cdot u \;.\]
\end{lemma}

\begin{proof}
By linearity and continuity, it suffices to show that $u\cdot\iota_f(a)=\iota_f(a)=\iota_f(a)\cdot u$ for any $f \in L$ and $a \in A_f$\/. Direct computation yields
\[ \begin{aligned}
u\cdot\iota_f(a)
& =\sum_{e\in L}\lm_e e\cdot\iota_f(a) \\
& = \sum_{e \in L} \lm_e \iota_{ef} \phi_{ef,f}(a) 
& = \sum_{h \in fL} \left(\sum_{e \st ef= h} \lm_e \right) \iota_h \phi_{h,f}(a)
\end{aligned} \]
(manipulations with sums are justified since $\sum_e \lm_e$ is an absolutely convergent series).

On the other hand, since $u$ is an identity for $\lp{1}(F)$ we have
\[ f = u\cdot f = \sum_{e \in L} \lm_e ef   = \sum_{h \in fL} \left(\sum_{e \st ef= h} \lm_e \right) h \]
and by comparing coefficients of $h$ on both sides of this equation we deduce that
\[ \sum_{e \st ef= h} \lm_e = \left\{ \begin{aligned} 0 & \quad\text{ if $h\neq f$} \\ 1 & \quad\text{ if $h=f$} \end{aligned} \right. \]
Therefore $u\cdot\iota_f(a) =\iota_f \phi_{f,f}(a) = \iota_f(a)$\/. The proof that $\iota_f(a)\cdot u=\iota_f(a)$ is identical save for switching left and right multiplication, and we omit the details.
\end{proof}

\subsection*{Multilinear extensions}
Later on, we shall need to check whether certain functions defined on ${\SA\times\ldots\times\SA}^n$ are multilinear. In the cases we need it is fairly obvious whether this is the case, but to be precise we include the following lemma. First we set up some terminology.
\begin{defn}\label{dfn:multilin_extn}
Let $\Ind$ be an indexing set and let $(X_i)_{i\in\Ind}$ be a family of Banach spaces. Let $E$ be a Banach space: if $\wtild{T}$ is a function $\left(\coprod_{k\in\Ind}X_k\right)^n \to E$, an \dt{extension of $\wtild{T}$ to $\X$} is an $n$-linear function
\[ T: \overbrace{\X\times\ldots\times\X}^n \to E \]
such that $T\circ(\inc{k_1}, \ldots ,\inc{k_n}) = \wtild{T}\vert_{X_{k_1}\times\ldots\times X_{k_n}}$ for all $k_1, \ldots, k_n \in \Ind$\/.
\end{defn}

\begin{lemma}\label{l:multilin_extn}
Let $n \in \Nat$ and let $C> 0$; let $E$ be a Banach space and let
\[ \wtild{T}:\left(\coprod_{i \in\Ind} X_i\right)^n \to E \]
be an arbitrary function.
Then $\wtild{T}$ has an extension to $\X$ of norm $\leq C$, if and only if the following condition is satisfied:

\vspace{0.5em}
\noindent{$\mathbf (*)$} \qquad
\begin{tabular}{l}
 for every $n$-tuple $\left(k_1, \ldots, k_n\right) \in \Ind^n$, the restriction of $\wtild{T}$ to \\
 ${X_{k_1} \times\ldots\times X_{k_n}}$ is a bounded $n$-linear function with norm $\leq C$\/.
\end{tabular}
\end{lemma}
Clearly if such an extension exists it will be unique, by linearity and continuity of the extension.

\begin{proof}[Proof of Lemma \ref{l:multilin_extn}]
If $\wtild{T}$ extends to $T$ and $\vctR{k}=(k_1,\ldots, k_n)\in \Ind^n$, then for each $j\in\{1,\ldots, n\}$ and every $(n-1)$-tuple $(x_1,\ldots, x_{j-1}, x_{j+1}, \ldots, x_n)$ where $x_r \in X_{k_r}$ for all $r$, the function $X_{j(r)}\to E$ given by
\[ y \mapsto \wtild{T}(x_1,\ldots, x_{j-1}, y ,x_{j+1}, \ldots, x_n) = T\left(\inc{k_1}x_1, \ldots, \inc{k_j}y,\ldots, \inc{k_n}x_n\right) \] is bounded linear with norm $\leq \norm{T}$\/. Hence $(*)$ holds.

Conversely, if $(*)$ holds we define the putative extension $T$ as follows. If $u_1, \ldots, u_n \in \X$ then each $u_r$ has a unique representation as an absolutely convergent sum
\[ u_r = \bigoplus_{k\in \Ind}\inc{k} x^r_k \qquad(x^r_k\in X_k\;\text{ for all $k\in \Ind$}) \]
and we make the definition
\[ T(u_1,\ldots, u_n) \defeq \sum_{\vctR{k}\in \Ind^n} \wtild{T}(x^1_{k_1},\ldots,x^n_{k_n})  \]
where the sum on the right-hand side is absolutely convergent and thus well-defined; moreover, we have by condition $(*)$ the following estimate:
\[ \sum_{\vctR{k}\in \Ind^n} \norm{\wtild{T}(x^1_{k_1},\ldots,x^n_{k_n})} \leq C \sum_{\vctR{k}\in \Ind^n} \norm{x^1_{k_1}}\dotsb \norm{x^n_{k_n}} = C \norm{u_1}\dotsb \norm{u_n} \;.\]

It remains only to show that $T$ is indeed multilinear. Let $\lm \in\Cplx$ and $r\in\{1,\ldots, n\}$; let $v\in \X$ have an expansion of the form $v=\sum_{k\in \Ind}\inc{k}y_k$ where $\norm{v}=\sum_{k\in \Ind}\norm{y_k}<\infty$; then
\[\begin{aligned}
& T(u_1,\ldots,\lm u_r+v, \ldots, u_n) \\
 & = \sum_{\vctR{k}\in \Ind^n} \wtild{T}(x^1_{k_1},\ldots,\lm x^r_{k_r}+y_{k_r},\ldots, x^n{k_n})  & \quad\text{(by definition of $T$)} \\
 & = \sum_{\vctR{e}\in \Ind^n} \left[
      \begin{aligned} & \lm\wtild{T}(x^1_{k_1},\ldots,x^r_{k_r},\ldots,x^n_{k_n}) \\
      & + \wtild{T}(x^1_{k_1},\ldots, y_{k_r},\ldots, x^n_{k_n}) \\ \end{aligned}
  \right] & \quad\text{(by Condition $(*)$)} \\
 & = \left\{
      \begin{aligned} & \lm T(u_1,\ldots,u_r,\ldots,u_n) \\
      & + T(u_1,\ldots, u_{r-1}, v,\ldots, u_n) \\ \end{aligned}
      \right. & \quad\text{(definition of $T$)}\\
 \end{aligned} \]
 and thus $T$ is linear in the $r$th variable. Since $r$ was arbitrary, $T$ is multilinear as required.
\end{proof}

\end{section}

\begin{section}{Hochschild homology for $\SA[L,A]$}
Having introduced the Banach algebras we shall study in the paper, we pause to give some general background on the homology theory that we will use.

The Hochschild \emph{cohomo\-logy theory} of Banach algebras has been studied by many authors, and we refer the reader to \cite{Hel_HBTA} and \cite{BEJ_CIBA} for the relevant background. 
Hochschild homology has perhaps received less overt attention.
Let us therefore briefly repeat the relevant definitions, in order to set out the notational conventions which will be used in this article. 

\begin{defn}
Let $A$ be a Banach algebra (not necessarily unital) and let $M$ be a Banach $A$-bimodule. For $n \geq 0$ we define
\begin{align*}
\Ho{C}{n}(A,M) & \defeq M\ptp A^{\ptp n} \\
\end{align*}
where $\ptp$ denotes the projective tensor product of Banach spaces.

For $0\leq i\leq n+1$ the \dt{face maps} $\face[n]{i}:\Ho{C}{n+1}(A,M)\to\Ho{C}{n}(A,M)$ are the contractive linear maps given by
\[ \face[n]{i}(x\tp a_1\tp\ldots\tp a_{n+1}) 
 = \left\{\begin{aligned}
      xa_1\tp a_2\tp \ldots\tp a_{n+1} & \quad\text{ if $i=0$} \\
      x\tp a_1\tp\ldots\tp a_i a_{i+1}\tp \ldots\tp a_{n+1} & \quad\text{ if $1\leq i \leq n$} \\
    a_{n+1}x\tp a_1 \tp \ldots \tp a_n & \quad\text{ if $i=n+1$}
 \end{aligned} \right. \]
and the \dt{Hochschild boundary operator} $\bdy_n:\Ho{C}{n+1}(A,M)\to\Ho{C}{n}(A,M)$ is given by
\[ \bdy_n = \sum_{j=0}^{n+1} (-1)^j \face[n]{j} \;. \]
With these definitions, the Banach spaces $\Ho{C}{n}(A,M)$ assemble into a chain complex
\[ \ldots \lTo^{\bdy_{n-1}} \Ho{C}{n}(A,M) \lTo^{\bdy_n} \Ho{C}{n+1}(A,M) \lTo^{\bdy_{n+1}} \ldots \]
called the \dt{Hochschild chain complex} of $(A,M)$\/.
We let
\begin{flalign*}
\qquad \Ho{Z}{n}(A,M) & \defeq \Ker \bdy_{n-1} & \text{(the space of \dt{$n$-cycles})} \\
\Ho{B}{n}(A,M) & \defeq \Image \bdy_n & \text{(the space of \dt{$n$-boundaries})} \\
\Ho{H}{n}(A,M) & \defeq  \frac{\Ho{Z}{n}(A,M) }{\Ho{B}{n}(A,M) } & \text{(the \dt{$n$th Hochschild homology group}\index{Hochschild homology})} 
\end{flalign*}
\end{defn}

\begin{rem}
In the literature the spaces defined above are often referred to as the space of \emph{bounded} $n$-chains, $n$-cycles and $n$-boundaries: the resulting homology groups are then called the \dt{continuous} Hochschild homology groups of $(A,M)$\/. In this article we have chosen to omit these adjectives for sake of brevity; this should not lead to any ambiguity.
\end{rem}

Later, in Section \ref{s:finfree}, we shall also need the notions of \dt{$K$-relative} homology groups where $K$ is a unital Banach algebra.
Let $A$ be a unit-linked Banach $K$-algebra (see Definition \ref{dfn:K-BA}).
\begin{defn}\label{dfn:rel-ho}
Let $X$ be a Banach $A$-bimodule. The space of \dt{$K$-normalised} $n$-chains on $A$ with coefficients in $X$ is the Banach space $X\ptpR{K^e}\left(A^{\ptp_K n}\right)$, and is denoted by $\Ho[K]{C}{n}(A,X)$\/.

More explicitly, $\Ho[K]{C}{n}(A,X)$ is defined to be the quotient of $X\ptp A^{\ptp n}$ by the closed subspace ${\mc N}_n(K)$, where ${\mc N}_n$ is the closed linear span of all tensors of the form
\[ \begin{aligned}
    &  xc \tp a_1\tp\ldots\tp a_n - x \tp ca_1 \ldots\tp a_n & \\
 \text{ or } & x\tp a_1 \ldots \tp a_jc\tp a_{j+1}\tp\ldots\tp a_n - x\tp a_1\ldots\tp a_j\tp ca_{j+1}\tp\ldots\tp a_n & \quad(1\leq j \leq {n-1}) \\
 \text{ or } &  x\tp a_1\tp\ldots\tp a_n c - cx\tp a_1\tp\ldots\tp a_n & \\
\end{aligned} \]
where $x \in X$, $c \in K$ and $a_1, \ldots, a_n\in A$\/.
\end{defn}

By considering the action of face maps on each such tensor, one sees that \\ $\bdy_{n-1}({\mc N}_n(K))\subseteq{\mc N}_{n-1}(K)$, and so ${\mc N}_\blob(K)$ is a subcomplex of the Hochschild chain complex $\Ho{C}{*}(A,X)$\/. Hence by a standard diagram chase the quotient spaces $\Ho[K]{C}{*}(A,X)$ form a quotient complex of $\Ho{C}{*}(A,X)$, whose homology groups are the \dt{$K$-relative homology groups} of $A$ with coefficients in $X$\/.

The $K$-relative homology groups of $A$ may be easier to compute than the `full' homology groups, but it is not always clear how to relate the two families. Under certain circumstances the natural maps $\Ho{H}{*}\to\Ho[K]{H}{*}$ are isomorphisms; we shall only require the following well-known instance of this phenomenon.

\begin{defn}\label{dfn:contractible}
A Banach algebra $K$ is said to be \dt{contractible} if there exists $\Delta\in K\ptp K$ satisfying the following conditions:
\begin{enumerate}[(i)]
\item $x\cdot  \Delta = \Delta\cdot x$ for all $x\in K$;
\item $x\pi(\Delta)=x=\pi(\Delta)x$ for all $x \in K$
\end{enumerate}
where $\pi: K\ptp K \to K$ is the product map. Such a $\Delta$, if it exists, is called a \dt{diagonal for $K$}.
\end{defn}

Let $K$ be a unital, contractible algebra and $\B$ a unit-linked Banach $K$-algebra. As a special case of standard results, we can `normalise chains on $\B$ with respect to~$K$'. Later we will need a quantitative version of this fact, stated below.

\begin{thm}\label{t:cntrct-normalise}
Let $K$ be a finite-dimensional, contractible Banach algebra with diagonal $\Delta$\/. Then there is a sequence of positive reals $C_n>0$ such that the following holds: whenever $A$ is a unit-linked Banach $K$-algebra and $X$ a Banach $A$-bimodule, there exists a chain map $\al: \Ho{C}{*}(A,X) \to \Ho{C}{*}(A,X)$ with the following properties:
\begin{enumerate}[$(a)$]
\item each $\al_n$ is $K$-normalised, i.e.~factors through the quotient map $\Ho{C}{n}(A,X)\to\Ho[K]{C}{n}(A,X)$;
\item there exists a chain homotopy from $\sid$ to $\al$, given by bounded linear maps $t_n:\Ho{C}{n}(A,X)\to\Ho{C}{n+1}(A,X)$ such that
$\bdy_n t_n+ t_{n-1}\bdy_{n-1}=\sid_n-\al_n$
and $\norm{t_n}\leq C_n$ for all $n$\/.
\end{enumerate}
(In particular, for each $n$ the canonical map $\Ho{H}{n}(A,X) \to \Ho[K]{H}{n}(A,X)$ is an isomorphism of semi\-normed spaces.)
\end{thm}
We have stated this result in a simple version suitable for our future purposes, and not in its full generality. In particular it is given for chains and not for cochains as this is the setting in which we shall need it.

The crucial observation for us is that our splitting homotopy may be bounded independently of $A$ and $X$; while the observation is far from new, it is not easy to find an explicit reference in the literature. A detailed proof may be found in Appendix A of the author's thesis \cite{YC_PhD}.

\subsection*{Statement of the disintegration theorem}
The aim of this article is to determine some of the cohomo\-logy groups for convolution algebras of the form $\SA[L,A]$\/. To make things clearer we adopt the following notation.

\begin{notn*}
Let $(L,A) \in \SAcat$\/. We shall write
\[ \begin{aligned}
\Ho{C}{*}[L;A] & \text{ for } & \Ho{C}{*}(\SA[L,A],\SA[L,A]) \\
\Ho{Z}{*}[L;A] & \text{ for } & \Ho{Z}{*}(\SA[L,A],\SA[L,A]) \\
\Ho{B}{*}[L;A] & \text{ for } & \Ho{B}{*}(\SA[L,A],\SA[L,A]) \\
\Ho{H}{*}[L;A] & \text{ for } & \Ho{H}{*}(\SA[L,A],\SA[L,A]) \\
\end{aligned}\]
(This notation is meant to suggest that the simplicial chain complex of $\SA[L,A]$ may be defined more directly in terms of the pair $(L,A)$, without introducing the intermediate object $\SA[L,A]$\/.)
\end{notn*}

In light of the main result in \cite{YC_GMJ} that $\lp{1}(L)$ has trivial simplicial homology for any semi\-lattice $L$, it is natural to enquire if one can recover the homology groups $\Ho{H}{*}[L;A]$ from the family of homology groups $\left( \Ho{H}{*}(A_e,A_e) \right)_{e\in L}$\/. To formulate this more precisely we make the following definition.

\begin{defn}
For each $n$ we let $\Ho[\diag]{C}{n}[L;A]$ denote the subspace
\[ \lpsum{1}_{e \in L} \Ho{C}{n}(A_e, A_e) \quad\subseteq \quad\Ho{C}{n}[L;A] \]
It is clear that $\Ho[\diag]{C}{*}[L;A]$ is a subcomplex of $\Ho{C}{*}[L;A]$\/. We shall sometimes refer to the elements of $\Ho[\diag]{C}{*}[L;A]$ as \dt{$L$-diagonal chains} on $\SA[L,A]$\/.
\end{defn}

Now we are able to state the main result of this paper.
\begin{thm}[Disintegration theorem]\label{t:mainresult}
Let $(L,A)\in\SAcat$\/. Then the inclusion of $\Ho[\diag]{C}{*}[L;A]$ into $\Ho{C}{*}[L;A]$ induces an isomorphism of homology groups.
\end{thm}
The point is, of course, that the complex $\Ho[\diag]{C}{*}[L;A]$ depends only on the underlying set of $L$ and the particular algebras $(A_e)_{e\in L}$, and \emph{not} on the transition homo\-morphisms $A_e\to A_f$ which are used to construct $\SA[L,A]$\/.

\begin{rem}
The proofs of \cite[Thm 2.1]{BowDunc} and \cite[Thm 4.6]{GouPouWh} may be easily extended to recover the cases $n=1$ and $n=2$ of Theorem \ref{t:mainresult}. The novelty of our result lies not so much in the generalisation from Clifford semi\-group algebras to strong semi\-lattices of algebras, as in systematically solving the higher-degree cohomo\-logy problems.
\end{rem}

We shall deduce Theorem \ref{t:mainresult} from a stronger but more technical one (Theorem \ref{t:mainthm} below), whose statement requires some more definitions. It turns out that the subcomplex $\Ho[\diag]{C}{*}[L;A]$ is a chain summand in $\Ho{C}{*}[L;A]$, and we are able to choose a chain projection onto it with good properties.

\subsection*{Splitting off the diagonal part}
Let $(L,A)\in\SAcat$, and denote the transition homo\-morphisms by $\phi_{f,e}:A_e \to A_f$ (where $f\preceq e$ in $L$). For each $n$, let $\mu^{L,A}_n: \Ho{C}{n}[L;A] \to \Ho[\diag]{C}{n}[L;A]$ be the linear contraction defined by
\[ \boxed{\mu^{L,A}_n(
\inc{e_0}a_0 \tp \ldots\tp\inc{e_n} a_n
) = \inc{p}\phi_{p,e_0}a_0\tp\ldots\tp\inc{p}\phi_{p,e_n}a_n} \]
where $e_0, \ldots, e_n \in L$, $a_j \in A_{e_j}$ for all $j$ and $p \defeq e_0\dotsb e_n$\/. We shall occasionally drop the superscript and write $\mu_n$ when it is clear which $(L,A)$ we are working with.

For each $n$, $\mu^{L,A}_n$ is a contractive, linear projection of $\Ho{C}{n}[L;A]$ onto $\Ho[\diag]{C}{n}[L;A]$ (this is clear by computation on block elements).

\begin{lemma}
$\mu^{L,A}_*$ is a chain map.
\end{lemma}
\begin{proof}
It suffices to prove that $\mu^{L,A}_*$ commutes with each face map on the simplicial chain complex $\Ho{C}{*}[L;A]$\/.

We do this in full detail for the face map $\face{0}$\/. Let $n \geq 1$; let $e_0, \ldots, e_n \in L$; and let $a_j \in A_{e_j}$ for $j=0,1,\ldots, n$\/. We write $p$ for $e_0e_1\dotsb e_n$ and $f$ for $e_0e_1$; then
\[ \begin{aligned}
 &  \mu^{L,A}_{n-1}\face{0}\left(\inc{e_0}a_0\tp\ldots\tp\inc{e_n}a_n\right)  &\\
 & = \mu^{L,A}_{n-1}\left((\inc{e_0}a_0\cdot\inc{e_1}a_1)\tp\ldots\tp\inc{e_n}a_n\right) & \\
 & = \mu^{L,A}_{n-1}\left(\inc{f}[\phi_{f,e_0}(a_0)\phi_{f,e_1}(a_1)]\tp\ldots\tp\inc{e_n}a_n\right) &
 \\
 & = \inc{p}\left[\phi_{p,e_0}(a_0)\phi_{p,e_1}(a_1)\right]\tp\ldots\tp\inc{p}\phi_{p,e_n}a_n & \quad
    \begin{gathered}\text{(defn of $\mu^{L,A}_{n-1}$, since} \\ \text{ $fe_2\dotsb e_n=p$)}
    \end{gathered} \\
 & = (\inc{p}\phi_{p,e_0}a_0) \cdot (\inc{p}\phi_{p,e_1}a_1)\tp\ldots\tp \inc{p}\phi_{p,e_n}a_n  & \quad
    \begin{gathered}\text{(since $\inc{p}$ is a}\\
    	\text{homo\-morphism)}\end{gathered} \\
 & = \face{0}\left( \inc{p}\phi_{p,e_0}a_0\tp\ldots\tp\inc{e}\phi_{p,e_n}a_n\right) & \\
 & = \face{0}\mu^{L,A}_n\left(\inc{e_0}a_0\tp\ldots\tp\inc{e_n}a_n\right)  & \quad
    \text{(defn of $\mu^{L,A}_n$)}\\
\end{aligned} \]
By linearity and continuity, we deduce that $\mu^{L,A}_{n-1}\face{0}=\face{0}\mu^{L,A}_n$ for all $n \geq 1$\/.

An exactly similar calculation shows that $\mu^{L,A}_{n-1}\face{i}=\face{i}\mu^{L,A}_n$, for each $i=1,2,\ldots, n$, and thus $\mu^{L,A}_{n-1}\bdy_{n-1}=\bdy_{n-1}\mu^{L,A}_n$ for all $n \geq 1$, as required.
\end{proof}

\begin{thm}\label{t:mainthm}
Let $\pi_n\defeq \sid_n - \mu_n : \Ho{C}{n}[L;A]\to \Ho{C}{n}[L;A]$\/. Then the chain projection $\pi_*: \Ho{C}{*}[L;A]\to\Ho{C}{*}[L;A]$ is null-homotopic: that is, there exist bounded linear maps $\sig_n: \Ho{C}{n}[L;A]\to\Ho{C}{n+1}[L;A]$ for each $n \geq 0$ such that
\[ \bdy_n\sig_n+\sig_{n-1}\bdy_{n-1}=\pi_n \quad\quad(n\geq 1) \]
and $\bdy_0\sig_0=\pi_0=0$\/.
\end{thm}

The methods used to prove Theorem \ref{t:mainthm} are inspired by calculations in the papers \cite{BowDunc, GouPouWh} for $\lp{1}$-algebras of Clifford semi\-groups; they constitute a more complicated implementation of the strategy used in the author's paper \cite{YC_GMJ}.
More precisely, as in \cite{YC_GMJ} we combine two lines of attack: the naturality of $\pi_*$ with respect to certain `transfer maps' (to be defined below); and the fact that we can prove the theorem in the special case where $L$ is a \emph{finite, free} semi\-lattice. We then set up an inductive argument to construct a splitting homotopy for $\pi_*$\/, by recursively `transferring' known splitting formulas from the finite free case over to \emph{natural} splitting formulas for the general case.

\end{section}

\begin{section}{Transfer along semi\-lattice homo\-morphisms}\label{s:transfer}
In this section we look at a particular class of homo\-morphisms between convolution algebras of the form $\SA[L,A]$, and at the chain maps thus induced between Hochschild chain complexes. We shall see later that these `transfer maps' allow us to make precise certain `changes of variable' that are required for later proofs.

Let $(L,A) \in \SAcat$\/. Suppose we have a semi\-lattice $H$ and a semi\-group homo\-morphism $\al: H \to L$\/. Then we obtain a semi\-lattice diagram $(H,A\al)$, by taking
\[ (A\al)_e\defeq A_{\al(e)} \qquad(e\in H) \]
and, whenever $e\succeq f$ in $H$, \emph{defining} the transition homo\-morphism
\[ (A\al)_e {\rTo^{\phi^{A\al}_{f,e}}} (A\al)_f\]
to be
\[ A_{\al(e)} {\rTo^{\phi^A_{\al f, \al e}}} A_{\al(f)}\/.\]
We may thus form the convolution algebra $(H,A\al)$\/. There is a contractive linear map $\tau_\al:\SA[H,A\al]\to \SA[L,A]$, defined by
\[ \tau_\al(\inc[A\al]{e}b) \defeq \inc[A]{\al e}b\quad\quad(e \in H; b \in (A\al)_e).\]
(Strictly speaking, $\tau_\al$ is defined
by applying Lemma \ref{l:multilin_extn} to
 the function $\wtild{\tau}: \coprod_{e\in H} (A\al)_e \to \SA[L,A]$\/, where
$\wtild{\tau}(b) = \inc[A]{\al e}(b)$ for $b\in (A\al)_e$\/.)

\medskip\noindent{\bf Claim.} $\tau_\al$ is an algebra homo\-morphism. 

\begin{proof}
Since $\tau_\al$ is linear and continuous it suffices to show that
\begin{equation}\label{eq:transfer_HM}
 \tau_\al(\inc[A\al]{e}b\cdot\inc[A\al]{f}c) = \tau_\al(\inc[A\al]{e}b)\cdot \tau_\al(\inc[A\al]{f}c)
\tag{{$\dagger$}}
\end{equation}
for all $b, c \in \coprod_{e\in H} (A\al)_e$, where $e\defeq\Base(b)$ and $f\defeq\Base(c)$\/.

Since
\[ \inc[A\al]{e}b\cdot\inc[A\al]{f}c\equiv\inc[A\al]{ef}\left[\phi^{A\al}_{ef,e}(b)\phi^{A\al}_{ef,f}(c)\right] \]
the left-hand side of \eqref{eq:transfer_HM} is, by definition of $\tau_\al$,
 \[ \tau_\al\left(\inc[A\al]{ef}[\phi^{A\al}_{ef,e}(b)\phi^{A\al}_{ef,f}(c)]\right) = \inc[A]{\al(ef)}\left[\phi^A_{\al(ef),\al e}(b)\phi^A_{\al(ef),\al f}(c)\right] \,; \]
while the right-hand side of \eqref{eq:transfer_HM} is, by the definitions of $\tau_\al$ and of multiplication in $\SA[L,A]$,
\[  \inc[A]{\al e}(b)\cdot \inc[A]{\al f}(c) =  \inc[A]{\al(e)\al(f)}\left[\phi^A_{\al(e)\al(f),\al(e)}(b) \phi^A_{\al(e)\al(f),\al(f)}(c) \right]\;. \]
The desired identity \eqref{eq:transfer_HM} now follows, since $\al(ef)=\al(e)\al(f)$\/.
\end{proof}

\begin{defn}
For $n \geq 0$, we let $\Tran{\al}_n$ denote the contractive linear map ${\tau_\al}^{\ptp n+1}: \Ho{C}{n}[H,A\al] \to \Ho{C}{n}[L,A]$\/. We call $\Tran{\al}_n$ the \dt{transfer map along $\al$} in degree $n$\/.
\end{defn}
Note that since $\tau_\al$ is a homo\-morphism, $\Tran{\al}_*$ is a chain map.

\begin{notn*}
Later on, it will be useful to have the following shorthand: if $\vctR{a}=(a_0, a_1,\ldots,a_n)$ is an $(n+1)$-tuple of elements of $\coprod_{e \in L}A_e$, where $a_i \in A_{e(i)}$, say, then we write
\[ \Vinc[A]{e}(\vctR{a}) \defeq \inc[A]{e_0}(a_0) \tp \inc[A]{e_1}(a_1)\tp \ldots \tp \inc[A]{e_n}(a_n) \]
and write $\vctR{e}=\Base(\vctR{a})$\/.
\end{notn*}
The next two lemmas are merely a matter of interpreting the new notation.

\begin{lemma}\label{l:basechange}
Let $L$, $H$, $\al$\/, and $A$ be as above. Then
$\Tran{\al}\Vinc[A\al]{e}(\vctR{a}) = \Vinc[A]{\al e}(\vctR{a})$
for every $e_0, \ldots, e_n\in L$ and $a_j\in A_{\al(e_j)}$\/, $j=0,\ldots, n$\/.
\end{lemma}

\begin{lemma}[`Transfer is functorial']\label{l:transfer_is_functorial}
Let $F,H,L$ be semi\-lattices and let
\[ F \rTo^\beta H \rTo^\al L \]
be semi\-group homo\-morhpisms. Then whenever $(L,A)\in\SAcat$,
\[ \Tran{\al}\Tran{}\Vinc[A\al\beta]{\blob} = \Tran{\al\beta}\Vinc[A\al\beta]{\blob}=\Vinc[A]{\al\beta\blob} \quad. \]
\end{lemma}

\begin{rem}
There is a more general notion of transfer, where we not only allow ourselves to change the base semi\-lattice but also the algebra structure above.

Specifically: let $(L,A), (H,B) \in \SAcat$ and let $\al:H \to L$ be a homo\-morphism of semi\-lattices. A \dt{transfer morphism $(H,B)\to (L,A)$} is then given by a family $(\mu_h: B_h \to A_{\al h})_{h \in H}$ of contractive, unital algebra homo\-morphisms such that the diagram
\[ \begin{diagram}[tight,height=2.5em]
B_h & \rTo^{\mu_h} & A_{\al h}\\
\dTo^{\phi^B_{h,j}} & & \dTo_{\phi^A_{\al h,\al j}} \\
B_j & \rTo_{\mu_j} & A_{\al j} 
\end{diagram} \]
commutes for every $h \leq j$ in $H$\/.

While this definition of transfer is in some sense `the right one' from the category-theoretic viewpoint, it is not needed for what follows and we shall not pursue it further.
\end{rem}

We now recall the chain projections $\mu^{L,A}$ that were defined earlier. It will be vital for the technical arguments used later that for fixed $n$, $\mu_n^{L,A}$ depends on $(L,A)$ in a well-behaved way; and here `well-behaved' means `compatible with transfer' in the following sense.
\begin{propn}[$\mu$ is natural]
$\mu^{\blob}_n$ commutes with transfer in degree $n$\/. That is: given $(L,A)\in\SAcat$, and a semi\-group homo\-morphism $\al: H \to L$ where $H$ is a semilattice, we have a commuting diagram
\[ \begin{diagram}[tight,width=4.5em]
 \Ho{C}{n}[H;A\al] & \rTo^{\mu^{H,A\al}_n} & \Ho{C}{n}[H;A\al] \\
 \dTo^{\Tran{\al}_n} & & \dTo_{\Tran{\al}_n} \\
 \Ho{C}{n}[L;A] & \rTo_{\mu^{L,A}_n} & \Ho{C}{n}[L;A] \\
\end{diagram} \]
\end{propn}
\begin{proof}
This is just book-keeping. Recall that $\Tran{\al}_n$ is just alternative notation for $\tau_\al^{\ptp n+1}$\/. Recall also that whenever $y\preceq x$ in $H$\/, $\phi^{A\al}_{y,x} = \phi^A_{\al y,\al x}$ by definition.

To keep the notation simple we shall give the proof in the special case $n=2$\/: it should be clear from this how the general case works.

We shall show that $\Tran{\al}_2\mu^{H,A\al}_2$ and $\mu^{L,A}_2\Tran{\al}_2$ coincide on block elements of $\Ho{C}{2}[H;A\al]$\/: by linearity and continuity this will imply that they coincide on every element of $\Ho{C}{2}[H;A\al]$\/.
 Thus, let $e,f,g\in H$ and let $a\in (A\al)_e$, $b\in (A\al)_f$, $c\in (A\al)_g$\/. 
Let $p=efg$\/: then, since $\al(p)=\al(e)\al(f)\al(g)$\/,
\begin{align*}
\Tran{\al}_2\mu^{H,A\al}_2(\inc[A\al]{e}a\tp\inc[A\al]{f}b\tp\inc[A\al]{g}c)
 & = \Tran{\al}_2 \left( \inc[A\al]{p}\phi^{A\al}_{p,e}(a) \tp \inc[A\al]{p}\phi^{A\al}_{p,f}(b) \tp \inc[A\al]{p}\phi^{A\al}_{p,g}(c)\right) \\
 & = \tau_\al\inc[A\al]{p}\phi^{A\al}_{p,e}(a) \tp \tau_\al\inc[A\al]{p}\phi^{A\al}_{p,f}(b) \tp \tau_\al\inc[A\al]{p}\phi^{A\al}_{p,g}(c) \\
 & = \inc[A]{\al p}\phi^{A\al}_{p,e}(a) \tp \inc[A]{\al p}\phi^{A\al}_{p,f}(b) \tp \inc[A]{\al p}\phi^{A\al}_{p,g}(c) \\
 & = \inc[A]{\al p}\phi^A_{\al p, \al e}(a) \tp \inc[A]{\al p}\phi^A_{\al p, \al f}(b) \tp \inc[A]{\al p}\phi^A_{\al p, \al g}(c) \\ 
 & = \mu^{L,A}_2 \bigl( \inc[A]{\al e}(a) \tp \inc[A]{\al f}(b) \tp \inc[A]{\al g}(c) \bigr) \\
 & = \mu^{L,A}_2 \bigl( \tau_\al\inc[A]{e}(a) \tp \tau_\al\inc[A]{f}(b) \tp \tau_\al\inc[A]{g}(c) \bigr) \\
 & = \mu^{L,A}_2 \Tran{\al}_2 ( \inc[A\al]{e}a\tp\inc[A\al]{f}b\tp\inc[A\al]{g}c) \\
\end{align*}
as required.
\end{proof}

\end{section}

\begin{section}{Explicit contractions in the finite free case}\label{s:finfree}
As in \cite{YC_GMJ} a key role is played by the special case of a finite, free semi\-lattice. Recall (Example \ref{eg:ff-slatt}) that if $\Ind$ is a finite set then the free semi\-lattice generated by $\Ind$ is the set of all non-empty subsets of $\Ind$, with the product of two such subsets defined to be their union.

\begin{propn}\label{p:F-normalise2}
Let $\Ind$ be a finite set and $F$ the free semi\-lattice generated by $\Ind$; let $(F,B)\in\SAcat$\/. Let $\mu^{F,B}_*:\Ho{C}{*}[F;B] \to \Ho{C}{*}[F;B]$ be the chain map defined earlier. Then there exists a chain homotopy from $\mu^{F,B}$ to $\sid$: more precisely, there exist bounded linear maps $s^F_n:\Ho{C}{n}[F;B]\to\Ho{C}{n+1}[F;B]$ such that
\[ \bdy^{F,B}_n s^F_n+ s^F_{n-1}\bdy^{F,B}_{n-1}=\sid_n-\mu^{F,B}_n\]
and such that $\norm{s^F_n}$ is bounded above by a constant depending only on $F$ and~$n$\/.
\end{propn}

The rest of this section is given over to a somewhat indirect proof of the proposition above, and so may be skipped if the reader wishes to get straight to the core of the proof of Theorem \ref{t:mainthm}.

Our approach will be to observe that one can identify the subcomplex of $L$-diagonal chains with the quotient complex of $\lp{1}(L)$-normalised chains, and to exploit the fact that if $F$ is a finite free semi\-lattice then $\lp{1}(F)$ is a contractible Banach algebra, to which we may apply Theorem \ref{t:cntrct-normalise}.

\subsection*{$L$-diagonal representatives for $\lp{1}(L)$-relative homology}
We work in slightly greater generality than is strictly necessary, as this seems to clarify the ideas involved.

We recall from the discussion after Definition~\ref{dfn:K-BA} that whenever $(L,A)\in\SAcat$\/, the convolution algebra $\SA[L,A]$ has the structure of an $\lp{1}(L)$-algebra. Therefore we may consider chains and cochains which are normalised with respect to this action of $\lp{1}(L)$\/.
For sake of legibility we shall write $\Ho[L]{C}{n}$
for the space of $\lp{1}(L)$-normalised chains.


\begin{lemma}\label{l:mu_normalises}
Let $(L,A)\in\SAcat$\/. For any $\vctR{a}\in\left(\coprod_{e\in L}A_e\right)^{n+1}$, $\Vinc{e}\vctR{a}$ and $\mu_n\Vinc{e}\vctR{a}$ have the same image under the canonical quotient map $q:\Ho{C}{*}[L;A]\to \Ho[L]{C}{*}[L;A]$\/.
\end{lemma}

\begin{proof}
We shall show how the proof works for the case $n=2$\/: the general case proceeds in the same way, but has to be formulated as a rather unwieldy inductive argument.

If $x$, $y$ are elementary tensors in $\Ho{C}{2}[L;A]$, let us temporarily say that $x$ and $y$ are \emph{$L$-equivalent} (denoted $x \sim y$) if $x-y \in \ker(q)$\/.

Let $\vctR{a}=(a_0, a_1, a_2)$ be an arbitrary element of $\left(\coprod_{e \in L} A_e\right)^3$; let $e_j=\Base(a_j)$ for each $j$\/. Recall that for any $e, f\in L$ and $a \in A_e$
\[ e\cdot\inc{e}a = \inc{e}a = \left(\inc{e}a\right) \cdot e
\qquad\text{ and }\qquad
 f\cdot \inc{e} a = ef\cdot\inc{e} a =\left(\inc{e}a\right)\cdot f \]

Let $p_2=e_2$, $p_1=e_1e_2$ and $p=e_0e_1e_2$\/.
By the definition of $L$-equivalence of tensors we see that
\[ \begin{aligned}
\Vinc{e}\vctR{a} & = \inc{e_0}a_0\tp\inc{e_1}a_1\tp\inc{e_2}a_2 \\
& = \inc{e_0}a_0\tp \inc{e_1}a_1\tp  \left[ e_2\cdot \inc{e_2}(a_2)\right] \\
& \sim \inc{e_0}a_0\tp \left[(\inc{e_1}a_1)\cdot e_2\right]\tp \inc{e_2}(a_2) \\
& = \inc{e_0}a_0\tp \left[ p_1\cdot \inc{e_1}a_1\right] \tp \left[p_2\cdot\inc{e_2}(a_2)\right]
\end{aligned} \]
and repeating this argument we obtain
\[\begin{aligned}
\Vinc{e}\vctR{a} & \sim
\inc{e_0}a_0\tp \left[ p_1\cdot \inc{e_1}a_1\right] \tp \left[p_2\cdot\inc{e_2}(a_2)\right] \\
& \sim \left[p\cdot\inc{e_0}a_0\right]\tp\left[p_1\cdot\inc{e_1}a_1\right]\tp \left[p_2\cdot\inc{e_2}(a_2)\right] \\
\end{aligned} \]
We may now run this argument in the opposite direction, `passing $p$ from left to right'. Since $pp_j=p$ for $j=0,1,2$\/, we get
\[ \begin{aligned}
&  \left[p\cdot\inc{e_0}a_0\right]\tp\left[p_1\cdot\inc{e_1}a_1\right]\tp \left[p_2\inc{e_2}(a_2)\right] \\
\sim &  \left[p\cdot\inc{e_0}a_0\right]\tp\left[p\cdot\inc{e_1}a_1\right]\tp \left[p_2\inc{e_2}(a_2)\right] \\
\sim & p\inc{e_0}a_0\tp p\inc{e_1}a_1\tp p\inc{e_2}a_2 \\
= &\inc{p}\phi_{p,e_0}a_0\tp\inc{p}\phi_{p,e_1}a_1\tp\inc{p}\phi_{p,e_2} a_2
\end{aligned} \]
and thus $\Vinc{e}\vctR{a} \sim \mu_n\Vinc{e}\vctR{a}$ as claimed.

\end{proof}

\subsection*{Constructing a homotopy for $\SA[F,B]$}

\begin{lemma}\label{l:F-gives-contractible}
Let $\Ind$ be a finite set and $F$ the free semi\-lattice generated by $\Ind$\/.
Then the Banach algebra $\lp{1}(F)$ is finite-dimensional and contractible.
\end{lemma}

This is well-known and we omit the proof, which can be found, for instance, as \cite[Lemma~5.5.1]{YC_PhD}\/; see also \cite[Example~1.6]{GHS_slatt}. We note that \cite{GHS_slatt} goes much further, providing an \emph{algorithm} for computing the diagonal of an \emph{arbitrary finite semi\-lattice}.

\begin{propn}\label{p:F-normalise1}
Let $F$ be as above and let $(F,B)\in\SAcat$; let $\B$ denote the convolution algebra $\SA[F,B]$; and let $X$ be a Banach $\B$-bimodule. Then there exists a chain map $\al: \Ho{C}{*}(\B,X) \to \Ho{C}{*}(\B,X)$ with the following properties:
\begin{enumerate}[$(a)$]
\item each $\al_n$ factors through the quotient map $q_n:\Ho{C}{n}(\B,X)\to \Ho[F]{C}{n}(\B,X)$;
\item there exists a chain homotopy from $\sid$ to $\al$, given by bounded linear maps $t_n:\Ho{C}{n}(\B,X)\to\Ho{C}{n+1}(\B,X)$ satisfying $\bdy_n t_n+ t_{n-1}\bdy_{n-1}=\sid_n-\al_n$ for all $n$;
\item the norm of each $t_n$ is bounded above by some constant depending only on $F$ and on $n$\/.
\end{enumerate}
\end{propn}

\begin{proof}
By Lemma \ref{l:F-gives-contractible} the algebra $\lp{1}(F)$ is contractible, so in particular has an identity element. By Lemma \ref{l:unitalSA}, the canonical inclusion of $\lp{1}(F)$ into $\SA[F,B]$ sends the identity element of $\lp{1}(F)$ to an identity element for $\SA[F,B]$, and the results $(a)$--$(c)$ now follow as a special case of Theorem \ref{t:cntrct-normalise}.
\end{proof}

In particular, on taking $X=\B$ we see that every simplicial cycle on $\SA[F,B]$ is homologous to an $F$-diagonal one. However, this is not quite enough:  we wish to show that every $n$-cycle $x\in \Ho{Z}{n}[F;B]$ is homologous to $\mu^{F,B}_nx$, so an extra step is needed.

This extra step can be done at a slightly more general level, and is given by the following trivial lemma. (It is stated for Hochschild homology but there is clearly a dual version for cohomo\-logy.)

\begin{lemma}[Combining normalising projections]\label{l:induced_null-htpy}
Let $A$ be a Banach algebra, $X$ a Banach $A$-bimodule. Suppose we have two chain maps $\al,\lm:\Ho{C}{*}(A,X) \to \Ho{C}{*}(A,X)$ such that
\begin{itemize}
\item $\al$ is chain-homotopic to the identity: i.e.~there exist bounded linear maps $t_n:\Ho{C}{n}(A,X)\to\Ho{C}{n+1}(A,X)$ such that $\bdy_n t_n+ t_{n-1}\bdy_{n-1}=\sid_n-\al_n$ for all $n$;
\item $\lm_*\al_*=\al_*$\/.
\end{itemize}
Then $\lm_*$ is chain-homotopic to the identity: more precisely, we have
\[ \bdy_n s_n+ s_{n-1}\bdy_{n-1}=\sid_n-\lm_n \]
where $\norm{s_n}\leq (1+\norm{\lm_{n+1}})\norm{t_n}$\/.
\end{lemma}
\begin{proof}
Let $n \geq 0$\/. We know that
\[ \bdy_n t_n+ t_{n-1}\bdy_{n-1}=\sid_n-\al_n\;;\] 
composing on the left with $\lm_n$ on both sides, and recalling that $\lm$ is a chain map, we have
\[ \bdy_n\lm_{n+1}t_n +\lm_nt_{n-1}\bdy_{n-1}=\lm_n-\lm_n\al_n=\lm_n-\al_n \]
Subtracting this equation from the previous one yields
\[ \bdy_n(\sid_{n+1}-\lm_{n+1})t_n + (\sid_n-\lm_n)t_{n-1}\bdy_{n-1} =\sid_n-\lm_n \;,\]
and on taking $s_n\defeq(\sid_{n+1}-\lm_{n+1})t_n$ we have proved the lemma.
\end{proof}

\begin{proof}[Proof of Proposition \ref{p:F-normalise2}]
By Proposition \ref{p:F-normalise1} there exists a chain map $\al:\Ho{C}{*}[F;B] \to \Ho{C}{*}[F,B]$ which factors through the quotient map $q: \Ho{C}{*}[F;B]\to\Ho[F]{C}{*}[F;B]$ and which is homotopic to the identity. Since $\al$ factors through~$q_n$\/, $\mu^{F,B}_n\al_n=\al_n$ for all~$n$, by Lemma~\ref{l:mu_normalises}. Hence Lemma~\ref{l:induced_null-htpy} applies.
\end{proof}

\begin{rem}
Note that for fixed $F$ the maps $s^{F,\blank}_n$ are natural in the second variable, in some sense \ldots but to make this precise, we need a more general notion of transfer.
\end{rem}

\end{section}

\begin{section}[The main splitting theorem]{Proving the main splitting theorem}
\subsection*{Formulating the inductive step}
The calculations of this section may seem rather messy. Figure \ref{fig:schematic} is meant to give some idea of the underlying strategy behind the details. 

\begin{figure}
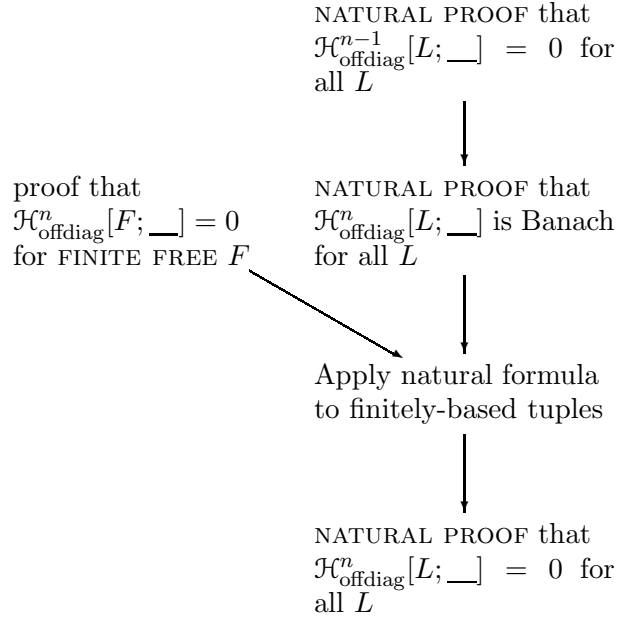

\begin{diagram}
& &  \parbox{4cm}{ {\sc natural proof} that \\ $\Co[{\rm offdiag}]{H}{n-1}[L;\blank]=0$ for all $L$}  \\
& &  \dTo \\
\parbox{4cm}{proof that \\ $\Co[{\rm offdiag}]{H}{n}[F;\blank]=0$ \\ for {\sc finite free} $F$} &  & \parbox{4cm}{{\sc natural proof} that \\ $\Co[{\rm offdiag}]{H}{n}[L;\blank]$ is Banach \\ for all $L$}  \\
& \rdTo &  \dTo  & \\
& & \parbox{4cm}{Apply natural formula \\ to finitely-based tuples} \\
& & \dTo \\
&  & \parbox{4cm}{ {\sc natural proof} that \\ $\Co[{\rm offdiag}]{H}{n}[L;\blank]=0$ for all $L$}   \\
\end{diagram}
\caption{The idea/motivation for the inductive step}
\label{fig:schematic}
\end{figure}

\begin{defn}
Let $j$ be a non-negative integer and let
\[ \left(\sig^{L,A}_j : \Ho{C}{j}[L;A] \longrightarrow \Ho{C}{j+1}[L;A] \right)_{(L,A) \in \SAcat} \]
be a family of bounded linear maps. We consider the following four conditions that $\sig^{\blob}$ may or may not satisfy.

\begin{itemize}
\item[\Cnd{R}] For each $(L,A)\in \SAcat$,
\[ \pi^{L,A}_{j+1}\sig^{L,A}_j = \sig^{L,A}_j \;.\]
\item[\Cnd{S}] For each $(L,A)\in\SAcat$
\[ \bdy^{L,A}_j \sig^{L,A}_j \bdy^{L,A}_j = \bdy^{L,A}_j\pi^{L,A}_{j+1} \;.\]
\item[\Cnd{T}] Whenever $(L,A)\in\SAcat$, $H$ is a semilattice and $\al: H \to L$ is a semi\-group homo\-morphism, then the diagram
\[ \begin{diagram}[tight,height=2.5em]
\Ho{C}{j}[H;A\al] & & \rTo^{\sig^{H,A\al}_j} & & \Ho{C}{j+1}[H;A\al] \\
\dTo^{\Tran{\al}_n} & &   & &\dTo_{\Tran{\al}_{n+1}} \\
\Ho{C}{j}[L;A] & & \rTo_{\sig^{L,A}_j} & & \Ho{C}{j+1}[L;A] \\
\end{diagram} \]
commutes.
\item[\Cnd{U}] There is a constant $K_j>0$ such that $\norm{\sig^{L,A}_j}\leq K_j$ for all $(L,A)\in\SAcat$\/.
\end{itemize}
\end{defn}
Here \Cnd{R} stands for {\sl range}\/ fixed by $\pi$; \Cnd{S} for {\sl semi\-lattice-normalised splitting}\/; \Cnd{T} for {\sl transferable}\/ map; and \Cnd{U} for {\sl uniformly bounded}\/.

\begin{rem}
Note that condition \Cnd{R} tells us that our splitting map $\sig_j$ should take `off-diagonal' cycles to `off-diagonal' cycles.
\end{rem}

\begin{lemma}[Base for induction]\label{l:induction_base}
For each $(L,A)\in\SAcat$ let $\sig^{L,A}_0: \Ho{C}{0}[L;A] \to \Ho{C}{1}[L;A]$ be the zero map. Then $\sig^\blob_0$ satisfies conditions \Cnd{R}--\Cnd{U}.
\end{lemma}
\begin{proof}
It is clear that taking $\sig^{L,A}_0=0$ for all $(L,A)\in\SAcat$ will satisfy conditions \Cnd{R}, \Cnd{T} and \Cnd{U}. To show that condition \Cnd{S} is satisfied we need only show that
\[ 0=d^{L,A}_0\pi^{L,A}_1\quad\quad\text{ for all $(L,A)\in\SAcat$\/.} \]

By linearity and continuity it suffices to verify this identity on block elements. So, let $e,f \in L$ and let $a\in A_e$, $b \in B_f$\/. Omitting superscripts for sake of clarity, we see that
\[ \mu_1(\inc{e}a\tp\inc{f}b) = \inc{ef}\phi_{ef,e}(a)\tp \inc{ef}\phi_{ef,f}(b)\equiv (f\cdot \inc{e}a)\tp (e\cdot\inc{f}b) \]
and so
\[ \begin{aligned}
\bdy_0\mu_1(\inc{e}a\tp\inc{f}b)
& = \bdy_0\left( (f\cdot \inc{e}a)\tp (e\cdot\inc{f}b)\right) \\
& = (\inc{f}b\cdot\inc{e}a) - (\inc{e}a\cdot\inc{f}b) 
& = \bdy_0(\inc{e}a\tp\inc{f}b)
\end{aligned} \]
which gives us $\bdy_0\pi_1(\inc{e}a\tp\inc{f}b)=\bdy_0(\inc{e}a\tp\inc{f}b)-\bdy_0\mu_1(\inc{e}a\tp\inc{f}b)=0$ as required.
\end{proof}

\begin{rem}It may help to think of the above proof as the `predual version' of the following statement: {\it every derivation $\SA[L,A] \to (\SA[L,A])'$ is automatically $L$-normalised}\/.\end{rem}

The following proposition is our inductive step, and is the heart of the main splitting theorem.
\begin{propn}\label{p:inductive_step}
Let $n \geq 1$, and suppose that there exists a family of linear maps $\sig^{L,A}_{n-1} : \Ho{C}{n-1}[L;A] \to \Ho{C}{n}[L;A]$ which satisfies conditions \Cnd{R}, \Cnd{S}, \Cnd{T} and \Cnd{U}.

Then there exists a family of linear maps $\sig^{L,A}_n : \Ho{C}{n}[L;A] \to \Ho{C}{n+1}[L;A]$ which satisfies conditions \Cnd{R}, \Cnd{T} and \Cnd{U}, and also satisfies
\begin{equation}\label{eq:condition_S'}
d^{L,A}_n\sig^{L,A}_n +\sig^{L,A}_{n-1}d^{L,A}_{n-1} = \pi_n^{L,A}\quad\quad\text{ for all $(L,A)\in\SAcat$}\;.
\tag{$\spadesuit$}
\end{equation}
\end{propn}

Let us first see how this proposition gives Theorem \ref{t:mainthm}.
\begin{proof}[Proof of Theorem \ref{t:mainthm}, using Proposition \ref{p:inductive_step}]
We shall prove a stronger statement, namely that there exists a family of bounded linear maps $\sig^{L,A}_n: \Ho{C}{n}[L;A]\to\Ho{C}{n+1}[L;A]$ for each $n \geq 0$ and $(L,A)\in\SAcat$, such that
\[ \bdy^{L,A}_n\sig^{L,A}_n+\sig^{L,A}_{n-1}\bdy^{L,A}_{n-1}=\pi^{L,A}_n \quad\quad(n\geq 1) \quad\quad\text{ and }\quad \bdy^{L,A}_0\sig^{L,A}_0=0 \]
for all $(L,A) \in \SAcat$, \emph{and} such that for each $n \in \Z_+$ the family $(\sig^{L,A}_n)_{(L,A)\in\SAcat}$ satisfies conditions \Cnd{R}, \Cnd{T} and \Cnd{U}.

The proof is by strong induction on $n$\/.

We have seen (Lemma \ref{l:induction_base}) that taking $\sig^\blob_0=0$ gives a family $(\sig^{L,A}_0)_{(L,A)\in\SAcat}$ which satisfies conditions \Cnd{R}, \Cnd{S}, \Cnd{T} and \Cnd{U}, and which satisfies
\[ \bdy^{L,A}_0\sig^{L,A}_0 = 0 \;.\]
Therefore by Proposition \ref{p:inductive_step} there exists a family $(\sig^{L,A}_1)_{(L,A)\in\SAcat}$ which satisfies conditions \Cnd{R}, \Cnd{T} and \Cnd{U}, and such that
\[ \bdy^{L,A}_1 \sig^{L,A}_1+\sig^{L,A}_0\bdy^{L,A}_0 = \pi^{L,A}_1 \]
for all $(L,A)\in\SAcat$\/. Hence
\[ \bdy^{L,A}_1\sig^{L,A}_1\bdy^{L,A}_1 = \left( \bdy^{L,A}_1 \sig^{L,A}_1+\sig^{L,A}_0\bdy^{L,A}_0 \right)\bdy^{L,A}_1 =\pi^{L,A}_1 \bdy^{L,A}_1 = \bdy^{L,A}_1\pi^{L,A}_2 \]
so that $\sig^\blob_1$ also satisfies condition \Cnd{S}.

Now suppose that there exists $n \geq 2$ and families $\sig^\blob_{n-2}$, $\sig^\blob_{n-1}$ which both satisfy conditions \Cnd{R}, \Cnd{T} and \Cnd{U}, and which also satisfy
\[ \bdy^{L,A}_{n-1}\sig^{L,A}_{n-1}+\sig^{L,A}_{n-2}\bdy^{L,A}_{n-2}=\pi^{L,A}_{n-1} \]
for all $(L,A)\in\SAcat$\/. Then
\[ \bdy^{L,A}_{n-1}\sig^{L,A}_{n-1}\bdy^{L,A}_{n-1} = \left( \bdy^{L,A}_{n-1} \sig^{L,A}_{n-1}+\sig^{L,A}_{n-2}\bdy^{L,A}_{n-2} \right)\bdy^{L,A}_{n-1} =\pi^{L,A}_{n-1} \bdy^{L,A}_{n-1} = \bdy^{L,A}_{n-1}\pi^{L,A}_n \]
and so $\sig^\blob_{n-1}$ satisfies condition \Cnd{S}. Hence by Proposition \ref{p:inductive_step} there exists a family $(\sig^{L,A}_n)_{(L,A)\in\SAcat}$ which satisfies conditions \Cnd{R}, \Cnd{T} and \Cnd{U} and which also satisfies
\[ \bdy^{L,A}_n\sig^{L,A}_n+\sig^{L,A}_{n-1}\bdy^{L,A}_{n-1}=\pi^{L,A}_n\quad. \]
This completes the induction, and hence Theorem \ref{t:mainthm} is proved \emph{assuming that Proposition \ref{p:inductive_step} is valid.}
\end{proof}

The proof of Proposition \ref{p:inductive_step} occupies the rest of this section.

\subsection*{Proof of the inductive step (Proposition \ref{p:inductive_step})}
\emph{Throughout the rest of this section} we fix $n \geq 1$ and let $F$ denote the finite semi\-lattice freely generated by $n+1$ idempotents $f_0, \ldots, f_n$\/.

Let $(F,T) \in \SAcat$\/. By Proposition \ref{p:F-normalise2} there exist bounded linear maps $s^T_n$ and $s^T_{n-1}$ such that:
\begin{equation}\label{eq:first_step}
\boxed{ \bdy^{F,T}_n s^T_n =\pi^{F,T}_n - s^T_{n-1} \bdy^{F,T}_{n-1} }
\tag{$\diamondsuit$}
\end{equation}
Recall also that $\norm{s^T_n}\leq C_n$ for some constant $C_n$ that is independent of $T$\/.

The formula \eqref{eq:first_step} says that every $n$-chain $z$ in the image of $\pi^{F,T}_n$ can be written as the sum of an $n$-boundary and a chain which is a linear image of $\bdy^{F,T}_{n-1}z$\/. The next step is aimed, roughly speaking, at improving this decomposition so that $z$ is the sum of an $n$-boundary and a chain which is the linear image of $\bdy^{F,T}_{n-1}z$ \emph{under a natural map} (Equation \eqref{eq:formal_identity} below).

Since $\pi^{F,T}$ is a chain projection, \eqref{eq:first_step} implies that
\begin{equation}\label{eq:second_step}
\bdy^{F,T}_n \pi^{F,T}_{n+1}s^T_n = \pi^{F,T}_n\bdy^{F,T}_n s^T_n =\pi^{F,T}_n - \pi^{F,T}_n s^T_{n-1} \bdy^{F,T}_{n-1} \;; 
\tag{$\dagger\dagger$}
\end{equation}
but, since condition \Cnd{S} holds for $\sig^{\blob}_{n-1}$, we also know that
\[ \bdy^{F,T}_{n-1}\left(\pi^{F,T}_n-\sig^{F,T}_{n-1}\bdy^{F,T}_{n-1} \right)=0 \;.\]
Hence post-multiplying on both sides of \eqref{eq:second_step} by $\pi^{F,T}_n-\sig^{F,T}_{n-1}\bdy^{F,T}_{n-1}$ gives us
 \begin{flalign*}
\qquad & \bdy^{F,T}_n  \pi^{F,T}_{n+1}s^T_n \left(\pi^{F,T}_n-\sig^{F,T}_{n-1}\bdy^{F,T}_{n-1} \right) & \\
& = \pi^{F,T}_n \left(\pi^{F,T}_n-\sig^{F,T}_{n-1}\bdy^{F,T}_{n-1} \right) & \\
& = \pi^{F,T}_n - \pi^{F,T}_n\sig^{F,T}_{n-1}\bdy^{F,T}_{n-1} & \quad\text{(since $\pi_n$ is a projection)} \\
& = \pi^{F,T}_n - \sig^{F,T}_{n-1}\bdy^{F,T}_{n-1} & \quad\text{(since \Cnd{R} holds for $\sig^\blob_{n-1}$)} 
\end{flalign*}
Before proceeding we  
introduce some auxiliary notation to ease the congestion of indices: let $\psi_T: \Ho{C}{n}[F;T] \to \Ho{C}{n}[F;T]$ denote the \emph{bounded linear} map
\[ \psi_T \defeq \pi^{F,T}_{n+1}s_n^T \left(\pi^{F,T}_n-\sig^{F,T}_{n-1}\bdy^{F,T}_{n-1} \right)  \]
This allows us to rewrite the rather cumbersome formula above more concisely, as
\begin{equation}\label{eq:formal_identity}
\boxed{ \bdy^{F,T}_n \psi_T=  \pi^{F,T}_n - \sig^{F,T}_{n-1}\bdy^{F,T}_{n-1}  }
\tag{$\heartsuit$}
\end{equation}
Observe for later reference that \label{eq:psi_is_normalised} $\pi^{F,T}_{n+1}\psi_T = \psi_T$\/. Also, since condition \Cnd{U} is assumed to hold, there is a constant $K_{n-1}$ such that $\norm{\sig^{F,T}_{n-1}}\leq K_{n-1}$\/. Therefore
\[ \norm{\psi_T} \leq \norm{s^T_n} (1+ K_{n-1}\norm{\bdy^{F,T}_{n-1} }) \leq C_n(1+(n+1)K_{n-1}) \]
where the right-hand side depends only on $n$\/.

\medskip
Now let $(L,A)\in\SAcat$\/. Let $x_0, \ldots, x_n \in L$, and observe that there is a unique, well-defined homo\-morphism of semi\-lattices $F \to L$ which sends $f_j \mapsto x_j$ for $j=0, \ldots, n$\/. We denote this homo\-morphism by $\ev{x}:F \to L$ --- the notation is meant to suggest that we think of $\ev{x}$ as `evaluation' of the `free variables' $f_j$ at particular values $x_j$ --- and write $\Ax$ for $A_\blob\circ\ev{x}: F \to \BAlg_1$\/.

Since we have a homo\-morphism $\ev{x}:F \to L$, we are in a position to bring in the transfer maps from Section \ref{s:transfer}. Consider the semi\-lattice diagram in $\BAlg_1$ given by $(F,\Ax)$, and the associated transfer chain map
\[ \Ho{C}{*}[F;\Ax] \rTo^{\Tran{\ev{x}}_* } \Ho{C}{*}[L;A] \]

Now suppose that for $j=0,1 \ldots, n$ we have $a_j\in (\Ax)_{f_j}=A_{x_j}$, so that
$\Vinc[A]{x}$ is an element of $\Ho{C}{n}[L;A]$\/. We shall define $\sig_n^{L,A}$ on elements of this form and extend using Lemma \ref{l:multilin_extn}.

Intuitively, Equation \eqref{eq:formal_identity} is a kind of formal identity in the unknowns $f_0, \ldots, f_n$ and $a_0, \ldots, a_n$, subject to the $f_j$ being commuting idempotents; therefore equality must be preserved when we `evaluate' each $f_j$ at $x_j$\/. With this observation in mind, we define a function $\wtild{\sig}^{L,A}_n: \left(\coprod_{x \in L} A_x \right)^{n+1} \to \Ho{C}{n+1}[L;A]$ as follows: given $\vctR{a} \in \left( \coprod_{x \in L}A_x\right)^{n+1}$, let
$\vctR{e}=\Base(\vctR{a})$ and set
\smallskip
\[ \boxed{ \wtild{\sig}^{L,A}_n(\vctR{a})\defeq \Tran{\ev{e}}_{n+1}  \psi_{A\ev{e}}\Vinc[A\ev{e}]{f}(\vctR{a}) } \]
We claim that $\wtild{\sig}^{L,A}_n$ has a bounded linear extension to $\SA[L,A]$\/. To see this, recall from Lemma \ref{l:multilin_extn} that $\wtild{\sig}^{L,A}_n$ has a linear extension of norm $\leq K$ if and only if, for every $\vctR{x}=(x_0,\ldots, x_n) \in L^{n+1}$, the restriction of $\wtild{\sig}^{L,A}_n$ to {$A_{x_0}\times\ldots\times A_{x_n}$ } is multilinear with norm $\leq K$\/. But by the way we have defined $\wtild{\sig}^{L,A}_n$, for fixed $\vctR{x}=(x_0,\ldots, x_n) \in L^{n+1}$
this restriction coincides with
$\Tran{\ev{x}}_{n+1}\psi_{\Ax}\Vinc[\Ax]{f}$\/, which is
 clearly multilinear since $\Vinc[\Ax]{f}$ is multilinear and $\Tran{\ev{x}}_{n+1}  \psi_{\Ax}$ is linear. Moreover, for any $(a_0, \ldots, a_n) $ we have
\[ \Norm{ \wtild{\sig}^{L,A}_n\Vinc[\Ax]{f} (a_0,\ldots, a_n)} \leq \norm{\psi_{\Ax} }\,\norm{a_0}\ldots\norm{a_n} \]
and we saw earlier that $\norm{\psi_{\Ax}} \leq K_n$ for some constant $K_n$ that depends only on $n$\/. Hence by Lemma \ref{l:multilin_extn} $\wtild{\sig}^{L,A}_n$ has a bounded $(n+1)$-linear extension to $\SA[L,A]$, as claimed. This extension may be canonically identified with a bounded \emph{linear} map $\Ho{C}{n}[L;A]\to \Ho{C}{n+1}[L;A]$, which we denote by $\sig^{L,A}_n$\/.

At this point, we have constructed for each $(L,A)\in\SAcat$ a bounded linear map $\sig^{L,A}_n: \Ho{C}{n}[L;A]\to \Ho{C}{n+1}[L;A]$, which satisfies $\norm{\sig^{L,A}_n}\leq K_n$ for some constant $K_n$ independent of $L$ and $A$\/. So in particular our family $\sig^\blob_n$ satisfies condition \Cnd{U}.

Next, we shall prove that the family $\sig^\blob_n$ satisfies condition \Cnd{T}. To do this we must show that the diagram
\[ \begin{diagram}[tight,height=2.5em]
\Ho{C}{n}[H;A\al] & & \rTo^{\sig^{H,A\al}_n} & & \Ho{C}{n+1}[H;A\al] \\
\dTo^{\Tran{\al}_n} & &   & &\dTo_{\Tran{\al}_{n+1}} \\
\Ho{C}{n}[L;A] & & \rTo_{\sig^{L,A}_n} & & \Ho{C}{n+1}[L;A] \\
\end{diagram} \]
commutes whenever we have a semi\-lattice homo\-morphism $H \rTo^\al L$\/.

By linearity and continuity it suffices to check this on block elements of $\Ho{C}{n}[H;A\al]$\/. Let $y_0$,\ldots, $y_{n+1}\in H$: for each $j$, let $x_j$ denote $\al(y_j) \in L$, and let $a_j \in (A\al)_{y_j}=A_{x_j}$\/. Recall that by Lemma \ref{l:transfer_is_functorial}
\[ \Tran{\al}_{n+1}\Tran{\ev{y}}_{n+1} = \Tran{\ev{x}}_{n+1} \]
and 
\[ \Tran{\al}_n\Vinc[A\al]{y}(\vctR{a}) = \Vinc[A]{x}(\vctR{a}) \;.\]
Therefore,
\[ \Tran{\al}_{n+1}\sig^{H,A\al}_n\Vinc[A\al]{y}(\vctR{a}) = \Tran{\al}_{n+1}\Tran{\ev{y}}_{n+1}\psi_{A\al\vctR{y}}\Vinc[A\al\ev{y}]{f}(\vctR{a}) = \Tran{\al\ev{y}}_{n+1}\psi_{A\al\ev{y}}\Vinc[A\al\ev{y}]{f}(\vctR{a}) \]
and
\[ \sig^{L,A}_n\Tran{\al}_n\Vinc[A\al]{y}(\vctR{a}) = \sig^{L,A}_n\Vinc[A]{x}(\vctR{a})= \Tran{\ev{x}}_{n+1}\psi_{\Ax}\Vinc[\Ax]{f}(\vctR{a}) \]
Since $\al\ev{y}=\ev{x}$, the previous two equations combine to give
\[ \Tran{\al}_{n+1}\sig^{H,B}_n\Vinc[A\al]{y}(\vctR{a}) = \sig^{L,A}_n\Tran{\al}_n\Vinc[A\al]{y}(\vctR{a}) \]
as required.

It remains only to verify conditions \Cnd{R} and \eqref{eq:condition_S'}. Let $(L,A)\in~\SAcat$: we must show that $\pi^{L,A}_{n+1}\sig_n=\sig^{L,A}_n$ and $\bdy^{L,A}_n\sig^{L,A}_n=\pi^{L,A}_n - \sig^{L,A}_{n-1}\bdy^{L,A}$\/.

Since we now know that $\sig^{L,A}_n$ is a bounded multilinear map $\Ho{C}{n}[L;A] \to \Ho{C}{n+1}[L;A]$, it suffices by linearity and continuity to check both putative identities on block elements of the form $\Vinc[A]{x}\vctR{a}$\/.

This is now mere diagram-chasing, given the machinery set up earlier. The diagram shown in Figure~\ref{fig:guiding_diag} is worth keeping in mind.
\begin{figure}[t]
\[ \begin{diagram}
\Ho{C}{n-1}[F;\Ax] & & \LR{\bdy^{F,\Ax}_{n-1}}{\sig^{F,\Ax}_{n-1}} & & \Ho{C}{n}[F;\Ax] & & \lTo^{\bdy^{F,\Ax}_n} & \Ho{C}{n+1}[F;\Ax] \\
\dTo^{\Tran{\ev{x}}_{n-1}} & & & & \dTo_{\Tran{\ev{x}}_n} & & &  \\
\Ho{C}{n-1}[L;A] & & \LR{\bdy^{L,A}_{n-1}}{\sig^{L,A}_{n-1}} & & \Ho{C}{n}[L;A] & & \pile{\lTo^{\bdy^{L,A}_n} \\ \rDash} & \Ho{C}{n+1}[L;A] \\
\end{diagram}\]
\caption{Transferring from $F$ to $L$}
\label{fig:guiding_diag}
\end{figure}
Since
\begin{flalign*}
\quad  \pi^{L,A}_{n+1}\sig^{L,A}_n (\Vinc[A]{x}\vctR{a})  
 & =\pi^{L,A}_{n+1} \Tran{\ev{x}}_{n+1} \psi_{\Ax} \Vinc[\Ax]{f}(\vctR{a}) & \quad\text{(by definition)} \\
 & = \Tran{\ev{x}}_{n+1} \pi^{F,\Ax}  \psi_{\Ax}\Vinc[\Ax]{f}(\vctR{a})  & \quad\text{($\pi$ commutes with transfer)} \\
& = \Tran{\ev{x}}_{n+1}  \psi_{\Ax}\Vinc[\Ax]{f}(\vctR{a})  &  \quad\text{(see remark after \eqref{eq:formal_identity})} \\
& = \sig^{L,A}_n (\Vinc[A]{x}\vctR{a})
\end{flalign*}
condition \Cnd{R} is satisfied. Finally (and \emph{this} is where we make use of condition \Cnd{T} in our induction):
\begin{flalign*}
\quad & \bdy^{L,A}_n\sig^{L,A}_n (\Vinc[A]{x}\vctR{a}) & \\
 & =\bdy^{L,A}_n \Tran{\ev{x}}_{n+1}   \psi_{\Ax} \Vinc[\Ax]{f}(\vctR{a}) & \text{(by definition)} \\
 & = \Tran{\ev{x}}_n \bdy^{F,\Ax}_n   \psi_{\Ax}\Vinc[\Ax]{f}(\vctR{a}) & \text{(transfer is a chain map)} \\
 & = \Tran{\ev{x}}_n \left( \pi^{F,\Ax}_n - \sig^{F,\Ax}_{n-1}\bdy^{F,\Ax}_{n-1} \right)(\Vinc[\Ax]{f}\vctR{a})  & \text{(by the `formal identity' \eqref{eq:formal_identity})} \\
& = \left(\Tran{\ev{x}}_n\pi^{F,\Ax}_n - \Tran{\ev{x}}_n\sig^{F,\Ax}_{n-1}\bdy^{F,\Ax}_{n-1} \right)(\Vinc[\Ax]{f}\vctR{a})  & \\
& = \left(\Tran{\ev{x}}_n\pi^{F,\Ax}_n - \sig^{L,A}_{n-1}\Tran{\ev{x}}_{n-1}\bdy^{F,\Ax}_{n-1} \right)(\Vinc[\Ax]{f}\vctR{a})  & \text{(condition \Cnd{T} for $\sig^\blob_{n-1}$)} \\
& =  \left( \pi^{L,A}_n\Tran{\ev{x}}_n- \sig^{L,A}_{n-1}\Tran{\ev{x}}_{n-1}\bdy^{F,\Ax}_{n-1} \right)(\Vinc[\Ax]{f}\vctR{a}) & \text{($\pi$ commutes with transfer)} \\
& =  \left( \pi^{L,A}_n\Tran{\ev{x}}_n- \sig^{L,A}_{n-1}\bdy^{L,A}_{n-1} \Tran{\ev{x}}_n \right)(\Vinc[\Ax]{f}\vctR{a}) & \text{(transfer is a chain map)} \\
& = \left( \pi^{L,A}_n - \sig^{L,A}_{n-1}\bdy^{L,A}_{n-1}  \right) \Tran{\ev{x}}_n (\Vinc[\Ax]{f}\vctR{a}) &  \\
& = \left( \pi^{L,A}_n - \sig^{L,A}_{n-1}\bdy^{L,A}_{n-1} \right) (\Vinc[A]{x}\vctR{a})
\end{flalign*}
where the last step follows from
 the definition of~$\ev{x}$ and Lemma~\ref{l:basechange}\/.

Thus the family $(\sig^{L,A}_n)_{(L,A) \in \SAcat}$ satisfies Equation \eqref{eq:condition_S'},  and this concludes the proof of Proposition \ref{p:inductive_step}. In view of the previous section, this completes the proof of Theorem \ref{t:mainthm} and hence of Theorem \ref{t:mainresult}.


\end{section}

\begin{section}{Applications to $\lp{1}$-convolution algebras}\label{s:Clifford}
\subsection*{Clifford semi\-groups}
Bowling and Duncan observed \cite[Thm 2.1]{BowDunc} that the $\lp{1}$-algebra of \emph{any} Clifford semi\-group $S$ is weakly amenable. By using Theorem \ref{t:mainresult} we can extend their result to higher-degree cohomo\-logy if we make further assumptions about~$S$\/.

\begin{propn}\label{p:Cliff_simpl-triv}
Let $S=\coprod_{e\in L}G_e$ be a Clifford semi\-group over the semi\-lattice~$L$\/. Suppose that each $G_e$ is amenable. Then  $\Ho{H}{n}(\lp{1}(S), \lp{1}(S))=0$ and $\Co{H}{n}(\lp{1}(S),\lp{1}(S)')=0$ for all $n \geq 1$\/.
\end{propn}

\begin{proof}
If $G_e$ is amenable then $\lp{1}(G_e)$ is amenable with constant~$1$\/, so that the simplicial chain complex $\Ho{C}{*}(\lp{1}(G_e,\lp{1}(G_e))$ is weakly split in degrees $1$ and above with constants independent of~$e$\/. Since
\[ \Ho[\diag]{C}{n}(\lp{1}(S),\lp{1}(S)) = \lpsum{1}_{e \in L} \Ho{C}{n}(\lp{1}(G_e),\lp{1}(G_e)) \]
it follows that $\Ho[\diag]{H}{n}(\lp{1}(S),\lp{1}(S))=0$ for all $n \geq 1$\/. Hence by Theorem \ref{t:mainresult},
\[ \Ho{H}{n}(\lp{1}(S),\lp{1}(S)) = \Ho[\diag]{H}{n}(\lp{1}(S),\lp{1}(S))=0 \quad\quad\text{ for all $n \geq 1$} \]
as required. The statement about cohomo\-logy follows by a standard duality argument, see e.g.~\cite[Coroll.~1.3]{BEJ_CIBA}\/.
\end{proof}

\begin{rem}
Note that in the degenerate case where $G_e=\{\id\}$ for each $e$, we recover the main result of \cite{YC_GMJ}, namely the simplicial triviality of the convolution algebra of a semi\-lattice.
\end{rem}

In \cite{YC_GMJ} vanishing results for the simplicial homology of $\lp{1}$-semi\-lattice algebras were used to deduce vanishing results for cohomo\-logy with arbitrary symmetric coefficients. This procedure carries over to the case of commutative Clifford semi\-groups.

\begin{thm}\label{t:comm_Cliff}
Let $S$ be a \emph{commutative} Clifford semi\-group and let $X$ be any symmetric Banach $\lp{1}(S)$-bimodule. Then $\Co{H}{n}(\lp{1}(S),X)=0$ for all $n \geq 1$\/.
\end{thm}

\begin{proof}
We first note that $\Co{H}{n}(\lp{1}(S),X)\iso\Co{H}{n}(\fu{\lp{1}(S)},X_1)$ where $X_1$ has underlying Banach space $X$ and is the natural \emph{unital} bimodule induced from~$X$\/. (This is a special case of a general result on cohomo\-logy of unitisations: see \cite[Exercise III.4.10]{Hel_HBTA} or \cite[Ch.~1]{BEJ_CIBA} for details.)

Clearly $\fu{\lp{1}(S)}=\lp{1}(\fu{S})$ where $\fu{S}$ is itself a (unital) commutative Clifford semi\-group. By Proposition \ref{p:Cliff_simpl-triv}, $\lp{1}(\fu{S})$ is simplicially trivial; hence the conditions of \cite[Thm 3.2]{YC_GMJ} are satisfied and we may use that theorem to deduce that $\Co{H}{n}(\fu{\lp{1}(S)},X_1) = 0$ for all $n \geq 1$\/, as required.
\end{proof}

\begin{rem}Note that the case $n=1$ follows from \cite[Thm 2.1]{BowDunc} and the well-known fact that weak amenability for commutative Banach algebras forces all bounded derivations with symmetric coefficients to vanish.
\end{rem}

What can be said for more general coefficients? Here matters are more delicate, and one difficulty seems to be the lack of a chain projection onto the subcomplex of $L$-diagonal chains (recall that in the case of \emph{simplicial} chains we could exploit the projection $\mu$ and its additional good properties).
For example, in \cite[Example 3.3]{BowDunc} the authors give a telling construction of a Clifford semi\-group $S$ with the following properties:
\begin{enumerate}
\item each constituent group of $S$ is amenable (in fact, can be taken to be the symmetric group on three objects);
\item there exists a non-inner bounded derivation $\lp{1}(S)\to \lp{1}(S)$\/.
\end{enumerate}
Thus in Theorem \ref{t:comm_Cliff} the condition of commutativity is essential.

\subsection*{Normal bands}

\begin{defn}
A semi\-group $S$ is said to be a \dt{band} if every element of $S$ is an idempotent. 
\end{defn}

A commutative band is nothing but a semi\-lattice, and one might hope that the techniques of \cite{YC_GMJ} extend to show simplicial triviality for the $\lp{1}$-algebras of bands.
While we have not been able to solve the general case, we can obtain positive results (in the sense of proving simplicial triviality) if we put further restrictions on our bands.

\begin{defn}
A band $R$ is said to be rectangular if $abc=ac$ for all $a,b,c \in R$\/. More generally, a \dt{normal band} is a band $S$ which satisfies the identity
\[ abca = acba \qquad (a,b,c \in S)\;. \]
\end{defn}

Normal bands were introduced in \cite{YamKim_nb}, where it was observed that they are precisely the class of bands which arise as strong semi\-lattices of rectangular bands. We may therefore apply the disintegration techniques developed in this paper.

\begin{lemma}\label{l:rect-band}
Let $R$ be a rectangular band. Then the chain complex\hfill\\ $\Ho{C}{*}(\lp{1}(R),\lp{1}(R))$ is split by bounded linear maps whose norm is independent of $R$\/.
\end{lemma} 

\begin{proof}
This essentially follows from the `$1$-biprojectivity' of $\lp{1}(R)$ -- that is, the existence of a contractive $\lp{1}(R)$-bimodule map $\rho: \lp{1}(R) \to \lp{1}(R)\ptp\lp{1}(R)$ which is right inverse to the product map. We can in fact give an explicit sequence of splitting maps as follows: fix $z \in R$, and define $s_n:\Ho{C}{n}(\lp{1}(R),\lp{1}(R)\to\Ho{C}{n+1}(\lp{1}(R),\lp{1}(R))$ to be the unique bounded linear map satisfying
\[ s_n( e_{x_0}\tp e_{x_1}\tp \ldots \tp e_{x_n}) \defeq e_{x_0z} \tp e_{zx_0} \tp e_{x_1} \tp\ldots \tp e_{x_n} \]
for every $x_0$, $x_1$, \ldots, $x_n$ in $R$\/. Clearly each $s_n$ is contractive, and it is easily verified that $\bdy_n s_n + s_{n-1}\bdy_{n-1}=\sid$ for all~$n$\/.
\end{proof}

\begin{propn}
Let $S$ be a normal band. Then the convolution algebra $\lp{1}(S)$ is simplicially trivial.
\end{propn}
\begin{proof}
Let ${\mc R}$ be the class of rectangular bands. By \cite[Propn IV.5.14]{How_fund-sgp} $S$ admits a strong semi\-lattice decomposition of type ${\mc R}$, say $S \iso \coprod_{e\in L} R_e$\/.

By Theorem \ref{t:mainresult} the inclusion
\begin{equation}\label{eq:FLINTOFF}
\lpsum{1}_{e\in L} \Ho{C}{*}(\lp{1}(R_e),\lp{1}(R_e)) \rSub \Ho{C}{*}(\lp{1}(S),\lp{1}(S))
\tag{$*$}
\end{equation}
induces isomorphism of homology groups; Lemma \ref{l:rect-band} implies that the complex on the left-hand side of \eqref{eq:FLINTOFF} is split exact in degrees $1$ and above.
\end{proof}

\begin{rem}
It seems likely that one could obtain a slightly more direct proof
of this result.
For, since
 one can construct `free normal bands' on any given generating set,
 the construction of splitting maps for arbitrary normal bands should
 follow from the finite free case, using
a recursive construction
as in~\cite{YC_GMJ}.
 The key point is to check that one has a `natural splitting map' in the sense of \cite{YC_GMJ} in degree $0$, in order to have a starting point for the induction.
\end{rem}
\end{section}

\subsection*{Acknowledgments}
This article is based on work from the author's PhD thesis \cite{YC_PhD} at the University of Newcastle upon Tyne, which was supported by an EPSRC grant. Particular thanks are due to M.~C.\ White, for suggesting that the techniques of \cite{YC_GMJ} should admit generalisations to the present context, and to Z.~A.\ Lykova and N.~Gr{\o}nb{\ae}k for a helpful critique of the original material.
The author also thanks N.~Spronk for providing a copy of the preprint \cite{GHS_slatt}.


\end{document}